\documentclass[11pt]{amsart}
\begin{document}

\newtheorem{theo}{Th\'eor\`eme}[section]
\newtheorem{prop}[theo]{Proposition}
\newtheorem{lemme}[theo]{Lemme}
\newtheorem{cor}[theo]{Corollaire}
\newtheorem{rem}[theo]{Remarque}
\numberwithin{equation}{section}
\leftskip -1cm
\rightskip -1cm

\def\cf{{\bf C}_L}

\def\ppl{{\bf P}}

\def\ssr{\'el\'ement semisimple r\'egulier}
\def\ssrs{\'el\'ements semisimples r\'eguliers}

\def\sgls{sous-groupe de Levi standard}
\def\sglss{sous-groupes de Levi standard}

\def\eci{r.e.c.i.}
\def\ecis{r.e.c.i.}

\def\ett{repr\'esentation essentiellement temp\'er\'ee}
\def\etts{repr\'esentations essentiellement temp\'er\'ees}

\def\rir{repr\'esentation irr\'eductible}
\def\rirs{repr\'esentations irr\'eductibles}

\def\rli{repr\'esentation lisse irr\'eductible}
\def\rlis{repr\'esentations lisses irr\'eductibles}
\def\g\'en\'eralis\'ee{g\'e\-n\'e\-ra\-li\-s\'ee}
\def\cusp{repr\'esentation cuspidale}
\def\cusps{repr\'esentations cuspidales}

\def\care{repr\'esentations de carr\'e int\'egrable}
\def\cares{repr\'esentations de carr\'e int\'egrable}

\def\temp{repr\'esentation temp\'er\'ee}
\def\temps{repr\'esentations temp\'er\'ees}

\def\aa{\mathcal A}
\def\z{\mathbb Z}
\def\cc{{\mathbb C}}
\def\G{\tilde{G}}
\def\g{{\mathbb G}}
\def\lra{\leftrightarrow}

\def\zzzf{\zeta_{FE}^l}
\def\zzz{\bar{\zeta}^l_{FE}}
\def\zzd{\bar{\zeta}_{D_FD_E}^l}
\def\zzzd{\zeta_{D_FD_E}^l}
\def\zzzz{\bar{\zeta}^{m''}_{FL}}
\def\zzzdl{\zeta_{D_FD_E}^l}
\def\zzdl{\bar{\zeta}_{D_FD_E}^l}

\def\lam{\lambda_{FE}^m}
\def\ll{\bar{\lambda}^m_{FE}}
\def\lll{\bar{\lambda}^m_{FE}}
\def\llam{\lambda_{EK}^m}
\def\llll{\bar{\lambda}_{EK}^m}
\def\llld{\bar{\lambda}_{D_FD_E}^m}
\def\lamd{\lambda_{D_FD_E}^m}
\def\lamdl{\lambda_{D_FD_E}^l}
\def\llldl{\bar{\lambda}_{D_FD_E}^m}

\def\vvf{v_{M_n(F)}}
\def\vvl{v_{M_n(L)}}

\def\T{{\bf{T}}_{F,l}}
\def\TT{\tilde{\bf{T}}_{F,l,A}}
\def\Tl{{\bf{T}}_{L,l}}
\def\TTl{\tilde{\bf{T}}_{L,l,\zzzdl(A)}}

\def\rep{re\-pr\'e\-sen\-ta\-ti\-on}

\title[Correspondance de Jacquet-Langlands \'etendue]{Correspondance
de  Jacquet-Langlands \'etendue \`a toutes les repr\'esentations}
\maketitle
\centerline{\footnotesize{par Alexandru Ioan
BADULESCU}{\footnote{Alexandru Ioan BADULESCU, Universit\'e de
Poitiers, UFR Sciences SP2MI, D\'epartement de Math\'ematiques,
T\'el\'eport 2, Boulevard Marie et Pierre Curie, BP 30179, 86962
FUTUROSCOPE CHASSENEUIL CEDEX\\
E-mail : badulesc\makeatletter @\makeatother
wallis.sp2mi.univ-poitiers.fr}}}
\ \\
\ \\

\def\ii{{\bf I}_{F,D}}

\def\ett{repr\'esentation essentiellement temp\'er\'ee}
\def\rir{repr\'esentation irr\'eductible}
\def\cusp{repr\'esentation cuspidale}
\def\care{repr\'esentation de carr\'e int\'egrable}
\def\temp{repr\'esentation temp\'er\'ee}
\def\lra{\leftrightarrow}

\def\temps{repr\'esentations temp\'er\'ees}
\def\cusps{repr\'esentations cuspidales}
\def\etts{repr\'esentations essentiellement temp\'er\'ees}
\def\cares{repr\'esentations de carr\'e int\'egrable}

\def\jl{{\bf JL}_r}
\def\rr{{\mathcal R}}
\def\jlh{{\bf JLH}}
\def\jla{{\bf JL}}
\def\ss{{\bf S}_{G,G'}}

\def\iii{{\bf i}}
\def\r{{\bf r}}
\def\sgls{sous-groupe de Levi standard}
\ \\
{\bf Abstract :} We give here a generalization of the local
Jacquet-Langlands correspondence to all the smooth irreducible
representations. This correspondence is characterized by the fact that
it respects the classical Jacquet-Langlands correspondence
and it commutes with the parabolic induction functor.
It has good behaviour with respect to the Jacquet functor and
the Aubert involution. Using this correspondence, we prove some
particular cases of the global Jacquet-Langlands correspondence and we
deduce that a class of representations of an inner form of $GL_n$
over a
$p$-adic field is unitarizable.\\


\tableofcontents

\section{Introduction}
\ \\

Le but principal de cet article est d'\'etudier plus a fond la
correspondance de Jacquet-Langlands, plus pr\'ecis\'ement chercher une
correspondance pour une classe plus large de repr\'esentations que
la classe des repr\'esentations de carr\'e int\'egrable. Non seulement
la question est naturelle, mais aussi elle intervient de fa\c{c}on
in\'eluctable quand on essaye de d\'emontrer une correspondance de
Jacquet-Langlands globale, entre les repr\'esentations automorphes de
$GL_n$ et les repr\'esentations automorphes
de ses formes int\'erieures sur un corps global. Na\"\i vement, on
attendrait
d'une telle correspondance qu'elle respecte localement, \`a chaque place,
la correspondance de Jacquet-Langlands locale. Ceci n'est en fait
jamais vrai. Prenons le cas le plus simple tra\^\i t\'e dans [Vi].
Soit $F$ un corps local non archim\'edien et $D$ une alg\`ebre \`a
division centrale et de dimension $n^2$ sur $F$. La correspondance de
Jacquet-Langlands locale r\'ealise une bijection entre l'ensemble des
classes d'\'equivalence des reper\'esentations de carr\'e int\'egrable
de $GL_n(F)$ et l'ensemble des classes d'\'equivalence des
repr\'esentations lisses irr\'eductibles unitaires de $D^*$. Mais, une
correspondance globale du type de celle \'etablie dans [Vi], mettra en
correspondance localement
une repr\'esentation lisse irr\'eductible unitaire de
$D^*$ soit avec la repr\'esentation de carr\'e int\'egrable $\pi$ de
$GL_n(F)$ qui lui correspond par Jacquet-Langlands, soit avec la duale
de $\pi$ au sens de [Ze]. Elle ne respecte donc
pas toujours localement la correspondance de Jacquet-Langlands.
Le cas g\'en\'eral est encore plus compliqu\'e,
une correspondance
globale compl\`ete faisant intervenir automatiquement et le spectre
cuspidal et le spectre r\'esiduel des deux groupes. Or, la correspondance
de Jacquet-Langlands locale classique ne tra\^\i te tout au plus que des
repr\'esentations temp\'er\'ees, tandis  que
les composantes locales des repr\'esentations automorphes
irr\'eductibles dans le spectre r\'esiduel de $GL_n$
ne sont pas temp\'er\'ees,
et souvent elles ne sont pas non plus duales des repr\'esentations de
carr\'e int\'egrable ([MW]). Donc, t\^ot ou tard se poseront les
questions,  si $G=GL_n(F)$, $F$ corps local non archim\'edien, et $G'$
est une forme int\'erieure de $G$,

1) en se donnant une repr\'esentation lisse
irr\'eductible (eventuellement unitaire) de $G$ qui
n'est pas temp\'er\'ee, y a-t-il une repr\'esentaion de 
$G'$ qui lui correspond ou non ? et

2) si oui, la repr\'esentation qui lui correspond est-elle
irr\'eductible ou non ? 

Notre article donne une r\'eponse assez satisfaisante
\`a la premi\`ere question, mais pour la deuxi\`eme question la
r\'eponse n'est pas satisfaisante, sauf dans le cas des
repr\'esentations duales (au sens de [Ze] ou [Au]) des
repr\'esentations temp\'er\'ees. Le lecteur remarquera que cette
classe de repr\'esentations est assez importante pour la
correspondance globale, vu la classification du spectre r\'esiduel de
$GL_n$ ([MW]). 

Soit $F$ un corps local non archim\'edien. Posons $G=GL_n(F)$ et soit
$G'$ une forme int\'erieure de $G$. Si $g\in G$ et $g'\in G'$ on dit
que $g$ et $g'$ {\it se correspondent} si $g$ et $g'$ sont semisimples
r\'eguliers et ont le m\^eme polyn\^ome caract\'eristique. Alors tout
\'el\'ement conjugu\'e \`a $g$ et tout \'el\'ement conjugu\'e \`a $g'$
se correspondent. Notons
$G_{G'}$ l'ensemble des \'el\'ements de $G$ qui correspondent \`a un
\'el\'ement de $G'$. Soient $f$ une fonction invariante par
conjugaison d\'efinie sur $G_{G'}$ et $f'$ une fonction invarainte par
conjugaison d\'efinie sur l'ensemble des \'el\'ements semisimples
r\'eguliers de $G'$. 
On dit que $f$ et $f'$  {\it se
correspondent} si elles sont \'egales sur des classes de conjugaion 
qui se
correspondent. Par abus, si $f$ est d\'efinie sur un sous-ensemble de
$G$
contenant $G_{G'}$, nous dirons qu'elle correspond \`a $f'$ si la
restriction de $f$ \`a $G_{G'}$ correspond \`a $f'$. Typiquement, ici,
les fonctions invariantes par conjugaison seront des combinaisons
lin\'eaires de caract\`eres de repr\'esentaions.

Alors la r\'eponse \`a la premi\`ere question pos\'ee plus
haut est la suivante : si $\pi$ est une repr\'esntation lisse
irr\'eductible de $G$, alors ou bien le caract\`ere 
de $\pi$ est nul
sur $G_{G'}$ ou bien sa restriction \`a $G_{G'}$ correspond \`a une 
combinaison lin\'eaire de caract\`eres
de repr\'esentations lisses
irr\'eductibles de $G'$ (le premier cas est en fait cas particulier du
deuxi\`eme, en consid\'erant la combinaison nulle). 
Les coefficents de la combinaison sont
entiers. Ceci se traduit par un morphisme de groupes entre les groupes de
Grothendieck $Grot(G)$ et $Grot(G')$
des repr\'esentations lisses de longueur finie de
$G$ et $G'$~:\\
\ \\
{\bf Th\'eor\`eme.} {\it Il existe un unique morphisme de groupes 
$${\bf LJ} :Grot(G)\to Grot(G')$$
tel que, pour tout $\pi\in Grot(G)$, le caract\`ere de $\pi$
correspond au caract\`ere de ${\bf LJ}(\pi)$. Ce morphisme est surjectif et
il prolonge la correspondance de Jacquet-Langlands au signe pr\`es.}\\

Nous pouvons \'egalement d\'ecrire le noyau de ${\bf LJ}$.
Ce qui est important aussi est de d\'efinir une section, un
morphisme en sens inverse. Ici il faut faire un choix.
Nous pouvons d\'efinir un morphisme injectif de groupes
$${\bf JL}:Grot(G')\to Grot(G),$$
prolongeant la correspondance de Jacquet-Langlands au signe pr\`es
et tel que, pour tout $\pi\in Grot(G')$, le caract\`ere de $\pi$
correspond au caract\`ere de ${\bf JL}(\pi)$.
Pour des \'enonc\'es plus pr\'ecis, voir le
th. \ref{jlr} et la prop. \ref{antiisomgrot}, section 3.1. Le choix
que nous avons fait pour d\'efinir ${\bf JL}$ semble tr\`es naturel,
comme nous l'expliquons plus bas, seulement ${\bf JL}$ ne transforme
pas une repr\'esentation irr\'eductible en une repr\'esentation
irr\'eductible. Le probl\`eme qui se pose, et auquel l'auteur de cet
article n'a pas de r\'eponse, est : existe-t-il une autre section de
${\bf LJ}$ qui transforme une repr\'esentation irr\'eductible en une
repr\'esentaion irr\'eductible ? Remarquons que les morphismes ${\bf
LJ}$ et ${\bf JL}$ tels que nous les avons d\'efinis, ensemble avec
l'\'etude de leur comportement vis-\`a-vis des foncteurs d'induction
et de restriction parabolique que nous faisons par la suite, semblent
\^etre importants pour r\'esoudre la question plus difficile de
l'irr\'eductibilit\'e. Par exemple, m\^eme si ${\bf LJ}$ ne respecte
pas la dualit\'e sur les deux groupes, il permet n\'eanmoins de
montrer que, si deux repr\'esentations se correspondent par ${\bf
LJ}$, les caract\`eres de
leurs duales se correspondent au sens d\'efini plus haut et au signe
pr\`es
(voir la
section \ref{inv}). Ce r\'esultat, appliqu\'e aux repr\'esentations duales
des repr\'esentations de carr\'e int\'egrable, joue un certain
r\^ole dans la compr\'ehension de la correspondance globale.

\`A partir des groupes de Grothendieck des repr\'esentations lisses de
longueur finie associ\'es aux groupes
lin\'eaires sur $F$, Zelevinski a construit une alg\`ebre de Hopf avec
involution ([Ze]). Le m\^eme type de construction peut \^etre faite
pour les groupes lin\'eaires sur une alg\`ebre \`a division $D$ de
dimension finie et
centrale sur $F$, les groupes \'etant cette fois des formes
int\'erieures du groupe lin\'eaire sur $F$ ([Ta]). Dans la sous-section 3.2,
nous montrons que les morphismes d\'efinis entre les groupes de
Grothendieck donnent lieu \`a des morphismes entre les alg\`ebres de
Hopf (la situation est un peu plus compliqu\'ee, voir le
th. \ref{jlh}). Ce morphisme entre les alg\`ebres de Hopf traduit le
comportement de la correspondance {\bf JL} que nous avons
d\'efinie vis-\`a-vis de l'induction parabolique et du foncteur de
Jacquet de restriction. C'est l'unique correspondance qui respecte les
caract\`eres, la
correspondance de Jacquet-Langlands classique, et commute avec
l'induction parabolique. Elle ne commute pas \`a proprement parlant
avec le foncteur de
Jacquet, mais les restrictions de deux repr\'esentations qui se
correpondent par $\bf JL$ ont des caract\`eres qui se correspondent.

Comme nous l'avons signal\'e, notre correspondance 
n'envoie pas (m\^eme au signe pr\`es)
les repr\'esentations irr\'eductibles sur
des repr\'esentations irr\'eductibles, et cela
m\^eme dans le cas le plus
simple o\`u $G'=GL_2(D)$, $D$ alg\`ebre \`a division, centrale et 
de dimension finie sur $F$ (voir le commentaire 2 au th. \ref{jlr}). 
Ceci constitue un probl\`eme pour 
l'\'etablissement
d'une correspondance de Jacquet-Langlands globale compl\`ete et il
faudra mieux comprendre le ph\'enom\`ene.  
D'autre part, la 
correspondance $\bf JL$ n'est pas uniquement d\'etermin\'ee par
la relation entre les
caract\`eres.
On peut donc se poser la question si, en renon\c{c}ant \`a
l'hypoth\`ese que la correspondance commute avec le foncteur
d'induction, on peut construire une autre correspondance qui envoie
des repr\'esentations irr\'eductibles sur des repr\'esentations
irr\'eductibles. Nous conjecturons l'existence d'une correspondance
qui v\'erifie la relation des caract\`eres et qui est d\'efinie par le
fait qu'elle
respecte les param\`etres de Langlands (sect. 3.4). 
Nous ne connaissons pas la
r\'eponse \`a cette question mais nous la v\'erifions dans le cas le plus
simple $G'=GL_2(D)$. 

Toutes ces constructions sont bas\'ees sur le fait que, si $G'$ est
une forme int\'erieure de $GL_n(F)$, $F$ corps local non
archim\'edien, 
alors l'induite
\`a $G'$ d'une repr\'esentation de carr\'e int\'egrable d'un
sous-groupe de Levi d'un sous-groupe parabolique de $G'$ est
irr\'eductible. Ce r\'esultat a \'et\'e prouv\'e par Deligne, Kazhdan
et Vign\'eras dans le cas o\`u la caract\'eristique de $F$ est
nulle. Suivant leur m\'ethode et utilisant la th\'eorie des corps
proches \`a la Kazhdan que nous avons developp\'ee dans [Ba2], nous
prouvons ce r\'esultat en caract\'eristique non nulle
(sect. 2). Ainsi, tous les r\'esultats sont valables en toute
caract\'eristique. 

Dans la section 4 nous montrons quelques cas particuliers de
correspondance de Jacquet-Langlands globale. Ici nous voyons
l'importance de la comprehension de la correspondance locale pour
d'autres repr\'esentations que les repr\'esentations de carr\'e
int\'egrable. Comme illustration du principe local-global-local, nous
obtenons comme corollaire de la correspondance globale le fait
 que les repr\'esentations duales des repr\'esentations de
carr\'e int\'egrable de $G'$ sont unitaires, ainsi que d'autres cas
d'unitarisabilit\'e, tous pr\'evus par les conjectures de [Ta]. 

Une partie importante de ce travail provient d'une th\`ese
\'ecrite sous la direction de Guy Henniart et je voudrais le
remercier ici pour son aide inestimable.
Je remercie \'egalement Colette Moeglin, Herv\'e Jacquet, Bertrand 
Lemaire et Anne-Marie Aubert qui ont lu partiellement
le manuscrit et lui ont 
apport\'e des corrections et suggestions. Merci \`a Andrei Moroianu
pour le contre-exemple de la section \ref{precision}.

\section{Irr\'eductibilit\'e des induites des representations de 
carr\'e int\'egrable}

Soient $F$ un corps local non archim\'edien de caract\'eristique
quelconque, $D$  une alg\`ebre \`a division centrale sur $F$, de
dimension finie $d^2$, $r$ un entier strictement positif et
$G'=GL_r(D)$.
Le but de cette section est de montrer le th\'eor\`eme
\ref{irred} plus bas. Le cas de caract\'eristique nulle a 
\'et\'e prouv\'e dans [DKV] (th\'eor\`eme B.2.d) et c'est l\`a
qu'appara\^\i t l'id\'ee d'utiliser le th\'eor\`eme de Paley-Wiener. 
On prouve ici le cas de caract\'eristique non nulle, en utilisant la 
m\'ethode des corps proches.
\begin{theo}\label{irred} Soit $P$ un sous-groupe parabolique
de $G'$ et soit $P=LU$ une d\'ecomposition de Levi de $P$. Soit $\pi$
une repr\'esentation de carr\'e int\'egrable de $L$. Alors
$ind_P^{G'}\pi$ est irr\'eductible. 
\end{theo}

On note $\Psi(G')$ l'ensemble des caract\`eres lisses non
ramifi\'es de $G'$. L'ensemble $\Psi(G')$ admet une structure naturelle de
vari\'et\'e alg\'ebrique. Si $\pi$ est une repr\'esentation admissible
de $G'$, on note $\Psi(G';\pi)$ l'ensemble de classes de
repr\'esentations du type $\psi\otimes\pi$, $\psi\in
\Psi(G')$ ($\Psi(G';\pi)$ a une structure de
vari\'et\'e alg\'ebrique isomorphe \`a un quotient de la vari\'et\'e
$\Psi(G')$). On note $Grot(G')$ le groupe de
Grothendieck des repr\'esentations lisses de longueur finie de
$G'$. \`A part l'ensemble des classes de repr\'esentations lisses
irr\'eductibles, $Grot(G')$ admet une autre $\z$-base
remarquable, voir paragraphes suivants. 

Appelons 
un sous-groupe parabolique de $G'$ {\it standard} s'il contient le
groupe des 
matrices triangulaires sup\'erieures, et un sous-groupe de Levi de $G'$ {\it
standard} s'il est un sous-groupe de Levi d'un sous-groupe parabolique
standard. Si $L$ est un \sgls\ de $G'$ alors il existe une
partition $A_1\coprod A_2\coprod ...\coprod A_k$ de l'ensemble $\{1,2\
...r\}$ o\`u $A_1=\{1,2\ ...r_1\},\ A_2=\{r_1+1,r_1+2\ ...r_1+r_2\}$
etc. telle que $L$ soit l'ensemble des matrices $M=(m_{ij})_{1\leq
i,j\leq r}\in G'$ telles que $m_{ij}$ est nul si $(i;j)\notin
\cup_{u=1}^kA_u\times A_u$. Nous
identifions alors $L$ au
produit $GL_{r_1}(D)\times GL_{r_2}(D)\times...\times
GL_{r_k}(D)$. Les notations du paragraphe pr\'ec\'edent s'appliquent aussi
\`a tout \sgls\ de $G'$.

Si $\nu\in \Psi(L)$, alors
$\nu=\prod_{u=1}^k\nu_u$, o\`u $\nu_u\in \Psi(GL_{r_u}(D))$ 
et s'\'ecrit donc $\nu_u=|\det|^{c_u}$, o\`u $c_u$ est un
nombre complexe. On dit que $\nu$ est strictement positif si
$u\mapsto Re(c_u)$ est une fonction strictement croissante.
 
Nous appelons {\it repr\'esentation de Langlands} une
repr\'esentation du 
type $ind_Q^{G'}\nu\otimes\tau$ o\`u $Q$ est un sous-groupe 
parabolique standard de $G'$, $\tau$ est une \temp\ du sous-groupe de
Levi de $Q$, et $\nu$ est un caract\`ere strictement positif. 
L'ensemble des repr\'esentations de Langlands forme 
une base 
de $Grot(G')$ ([DKV], A.4.f) que nous appellerons {\it la base de
Langlands}. Nous appellerons {\it repr\'esentation de Langlands
strictement induite} une repr\'esentation comme plus haut avec $Q\neq
G'$. On note $W(G')$ le
groupe form\'e par les matrices de permutation dans $G'$ qu'on
identifie aussi avec le groupe des permutations de l'ensemble $\{1,2\
...r\}$. Si $L$ est un \sgls\ de $G'$ et $g\in G'$, 
alors on pose $^gL=gLg^{-1}$
(m\^eme 
notation aussi pour les paraboliques) et
si $\pi$ est une repr\'esentation de $L$, on note $^g\pi$ la
repr\'esentation de $^gL$ d\'efinie par $^g\pi(x)=\pi(g^{-1}xg)$. On
utilise aussi les notations~: $L^g$ pour $^{g^{-1}}L$ et $\pi^g$ pour
$^{g^{-1}}\pi$. 

Une repr\'esentation de $G'$ est dite
{\it essentiellement
de carr\'e int\'egrable} si elle peut s'\'ecrire comme
produit tensoriel d'une repr\'esentation de carr\'e
int\'egrable et d'un caract\`ere de $G'$. 
Nous abr\'egeons ici
repr\'esentation(s)  essentiellement
de carr\'e int\'egrable  par {\it r.e.c.i.}. 
Nous en rencontrerons beaucoup par la suite. Une repr\'esentation est
dite {\it essentiellement temp\'er\'ee} si elle peut s'\'ecrire comme
produit tensoriel  d'une repr\'esentation
temp\'er\'ee et d'un caract\`ere  de $G'$.
\\
\ \\
{\bf D\'emonstration.} On peut supposer que $P$ est standard. Nous
montrons le th\'eor\`eme par r\'ecurrence sur 
l'entier strictement positif $k$ tel que $G'=GL_k(D)$ ($D$ est
fix\'ee). Pour $k=1$ le th\'eor\`eme est 
\'evident. Supposons que le th\'eor\`eme est v\'erifi\'e pour tout
$k<r$. Supposons maintenant par l'absurde qu'il existe un sous-groupe
parabolique propre $P_0$ de $G'=GL_r(D)$, une d\'ecomposition de Levi
$P_0=L_0U_0$ de $P_0$ et une repr\'esentation de carr\'e int\'egrable
$\pi_0$ de $L_0$ telle que l'induite de $P_0$ \`a $G'$ de la
repr\'esentation $\pi_0$ ne soit pas irr\'eductible. On sait qu'on a
dans $Grot(G')$~: 

$$ind_{P_0}^{G'}\pi_0=\sum_{i=1}^s a_i\tau_i$$
o\`u les $a_i$ sont des entiers strictement positifs et les $\tau_i$
sont des \temps\ de $G'$ non \'equivalentes, et aucune parmi les
$\tau_i$ n'est de carr\'e int\'egrable. 

\begin{lemme}\label{s=1} On a  $s=1$. 
\end{lemme}

{\bf D\'emonstration.} Supposons par l'absurde que  $s\geq 2$. On a
alors deux cas~:\\ 
\ \\
{\bf Premier cas~:} Il existe un $i\in \{2,3...s\}$ tel que
$\tau_i\notin \Psi(G';\tau_1)$. On peut supposer que $i=2$. Posons
$$\{\tau_1,\ \tau_{2},\ \tau_{3}\ ...\tau_{s}\}\cap \Psi(G';\tau_1)=
\{\tau_{i_1}=\tau_1,\ \tau_{i_2},\ \tau_{i_3}\ ...\tau_{i_p}\}$$
et 
$$\{\tau_1,\ \tau_{2},\ \tau_{3}\ ...\tau_{s}\}\cap\Psi(G';\tau_2)=
\{\tau_{j_1}=\tau_2,\ \tau_{j_2},\ \tau_{j_3}\
...\tau_{j_q}\}.$$ 
Soient
$$\alpha=\sum_{u=1}^qa_{j_u}$$
et
$$\beta=\sum_{v=1}^pa_{i_v}.$$
On consid\`ere alors la forme lin\'eaire $f$ sur $Grot(G')$ d\'efinie
sur la base de Langlands par~: 

1) $f(\rho)=\alpha$ si $\rho\in \Psi(G';\tau_1)$

2) $f(\rho)=-\beta$ si $\rho\in \Psi(G';\tau_2)$

3) $f(\rho)=0$ si $\rho$ est une repr\'esentation de la base de
   Langlands de $Grot(G')$ qui n'est pas dans l'un
   des deux cas pr\'ec\'edents.\\  
\ \\
{\bf Remarque.} Les conditions 1) et 2) ont \'et\'e choisies de
fa\c{c}on \`a ce que 1), 2) et 3) impliquent que $f(\chi\otimes
ind_{P_0}^{G'}\pi_0)=0$ pour tout $\chi\in \Psi(G')$. On a aussi
$f(\tau_1)=\alpha\neq 0$.\\ 
\ \\
{\bf Deuxi\`eme cas~:} Pour tout $i\in \{2,3...s\}$ il existe un
caract\`ere $\chi_i\in \Psi(G')$ tel qu'on ait
$\tau_i=\chi_i\otimes \tau_1$. Alors, par l'unicit\'e modulo
conjugaison de l'induite d'une \care\ qui contient une \temp\ donn\'ee,
on en d\'eduit que les $\chi_i$ v\'erifient~: 
$$(res_{P_0}^{G'}\chi_i)\otimes \pi_0 \simeq \pi_0^w$$
o\`u $w$ est un \'el\'ement du groupe de Weyl du tore maximal standard
de $L_0$ dans $G'$. Soit $K$ l'ensemle de tous les caract\`eres 
$\chi\in \Psi(G')$ qui v\'erifient une telle relation. $K$ est un 
sous-groupe de $\Psi(G')$.
Si  $\chi_1$ est le
caract\`ere trivial de $G'$, $\{\chi_1,\chi_2...\chi_s\}$ admet 
une structure de
groupe multiplicatif, quotient de $K$ par le stabilisateur de la classe
$\tau_1$.
En particulier, tous les $a_i$ sont \'egaux. 
Soit alors $f_0$ une fonction alg\'ebrique
d\'efinie sur la vari\'et\'e $\Psi(G';\tau_1)$ et qui v\'erifie~: 

- $f_0(\tau_1)=1$

- $f_0(\tau_i)=0$ pour $i\in\{2,3...s\}$.

Soit $\epsilon$ une racine primitive d'ordre $s$ de l'unit\'e. Posons~:
$$f_1(\tau)=\sum_{i=1}^s \epsilon^{i-1}f_0(\chi_i\otimes \tau)\ \ \ \
\ \ \ \forall \tau\in \Psi(G';\tau_1).$$ 
D\'efinissons cette fois la fonction $f$ sur la base de Langlands de
$Grot(G')$ de la fa\c{c}on suivante~: 

1') $f(\rho)=f_1(\rho)$ si $\rho\in \Psi(G';\tau_1)$

2') $f(\rho)=0$ si $\rho$ est une repr\'esentation de la base de
Langlands qui ne se trouve pas dans $\Psi(G';\tau_1)$.\\ 
\ \\
{\bf Remarque.} Dans ce cas aussi, $f(\chi\otimes
ind_{P_0}^{G'}\pi_0)=0$ pour tout $\chi\in \Psi(G')$, tandis que
$f(\tau_1)=1\neq 0$.\\ 

Nous allons montrer que la fonction $f$ v\'erifie toujours (qu'on soit
dans le premier ou le deuxi\`eme cas) les conditions de
Paley-Wiener. 
\begin{lemme}\label{fannule} L'application $f$
s'annule sur toute repr\'esentation induite \`a partir d'un
sous-groupe parabolique propre de $G'$.
\end{lemme}

{\bf D\'emonstration.} Se fait en deux \'etapes~:\\
\ \\
{\it \'Etape 1~: $f$ s'annule sur toute induite d'une repr\'esentation
essentiellement de carr\'e int\'egrable d'un sous-groupe de Levi
propre. (N'utilise pas l'hypoth\`ese de r\'ecurrence.)}\\ 

Soient $P$ un sous-groupe parabolique propre de $G'$, $P=LU$ une
d\'ecomposition de Levi de $P$ et $\sigma$ une \eci\ de $L$. On sait
(voir, par exemple, [Ba1], lemme 2.3) qu'on a deux possibilit\'es~:

- si $\sigma$ est un produit du type $\sigma=\psi\otimes \sigma_u$
  o\`u $\psi$ est la restriction \`a $L$ d'un caract\`ere $\Psi$ de
  $G'$ et $\sigma_u$ est de carr\'e int\'egrable,  alors
  $ind_L^{G'}\sigma=\Psi\otimes ind_L^{G'}\sigma_u$,  

- sinon $ind_L^{G'}\sigma$ est  une somme de repr\'esentations
  de Langlands strictement induites.  

Maintenant, si $\sigma$ est dans la seconde situation,
$f(ind_{P}^{G'}\sigma)=0 $ par la construction de $f$ (condition 3)
dans le premier cas et 2') dans le deuxi\`eme). Si $\sigma$ est dans
la premi\`ere situation, alors on a encore une fois deux
possibilit\'es~: 

- ou bien $L=L_0^g$ pour un $g\in G'$  et il existe un caract\`ere
  $\theta$ de $L$ qui est la restriction d'un caract\`ere $\Theta$ non
  ramifi\'e de $G'$ tel que $\sigma_u=\theta\otimes \pi_0^g$~; alors,
  dans $Grot(G')$, on a $ind_L^{G'}\sigma_u=\Theta \otimes
  ind_{L_0}^{G'}\pi_0$, ce qui implique  
$$f(ind_{P}^{G'}\sigma)=f(\Psi\otimes ind_L^{G'}\sigma_u)=f(\Psi
\Theta \otimes ind_{L_0}^{G'}\pi_0)=0$$ 
 (voir les remarques faites plus haut sur la construction de $f$ dans les
deux cas),  

- ou bien on ne se trouve pas dans cette situation~; alors les
s\'eries de composition de $\Theta\otimes ind_{P_0}^{G'}(\pi_0)$ et
$ind_P^{G'}\sigma_u$ sont disjointes pour tout $\Theta\in
\Psi(G')$. La suite de composition de $ind_P^{G'}\sigma_u$ est donc
form\'ee de repr\'esentations temp\'er\'ees, mais  ne contient aucun
\'el\'ement de $\cup_{i=1}^s\Psi(G';\tau_i)$. Cela implique que la
suite de composition de $ind_P^{G'}\sigma=\Psi\otimes
ind_P^{G'}\sigma_u$ est form\'ee de \etts\, mais ne contient aucun
\'el\'ement de $\cup_{i=1}^s\Psi(G';\tau_i)$. Donc, encore une fois,
$f(ind_{P}^{G'}\sigma)=0 $ par la condition 3) dans le premier cas et
2') dans le deuxi\`eme cas dans la construction de $f$.\\ 
\ \\ 
{\it \'Etape 2~:  Soit $Grot_{ind}(G')$ le sous-groupe de $Grot(G')$
engendr\'e par toutes les repr\'esentations induites \`a partir de
sous-groupes paraboliques propres de $G'$. Les repr\'esentations
induites des \ecis\ \`a partir 
de sous-groupes paraboliques propres de $G'$ forment une famille
g\'en\'eratrice pour  $Grot_{ind}(G')$. (Utilise l'hypoth\`ese de
r\'ecurrence.)}\\   

Remarquons que l'hypoth\`ese de r\'ecurrence implique que pour tout
sous-groupe parabolique propre $P$ de $G'$ qui a une d\'ecomposition
de Levi $P=LU$, toute repr\'esentation  temp\'er\'ee de $L$ est une
repr\'esentation induite d'une repr\'esentation de carr\'e
int\'egrable. Par cons\'equent, toute repr\'esentation
essentiellement temp\'er\'ee de $L$ est une repr\'esentation induite
d'une repr\'esentation essentiellement de carr\'e int\'egrable.  

Soit $P=LU$ un sous-groupe parabolique propre de $G'$. La remarque
plus haut vaut maintenant pour tous les sous-groupes paraboliques de
$L$, propres ou pas cette fois. Mais l'ensemble des induites des
\etts\ de tous les sous-groupes paraboliques (propres ou pas) de $L$
est une famille g\'en\'eratrice de $Grot(L)$ (car elle contient la
base de Langlands de $L$). L'hypoth\`ese de r\'ecurrence implique donc
que l'ensemble des induites des \ecis\ de tous les sous-groupes
paraboliques (propres ou pas) est aussi une famille g\'en\'eratrice de
$Grot(L)$. Cela prouve que les repr\'esentations induites de \ecis\
\`a partir de sous-groupes paraboliques propres de $G'$  engendrent
$Grot_{ind}(G')$.\qed
\ \\

Pour v\'erifier les
conditions de Paley-Wiener sur $f$, il suffit, gr\^ace au lemme
\ref{fannule}, 
de montrer que pour   
toute repr\'esentation irr\'eductible $\pi$ de $G'$, la restriction de
$f$ \`a $\Psi(G';\pi)$ est alg\'ebrique. Pour cela, on \'ecrit $\pi$
sur la base de Langlands dans $Grot(G')$. Il y a deux types de
repr\'esentations qui apparaissent dans cette \'ecriture~: des \etts\
de $G'$ et des repr\'esentations de Langlands induites strictes. Quand
on fait le produit tensoriel d'une  repr\'esentation de Langlands
induite stricte par un caract\`ere, on obtient toujours un \'el\'ement
de $Grot_{ind}(G')$. Or, on a montr\'e que $f$ s'annule sur
$Grot_{ind}(G')$. Donc, l'alg\'ebricit\'e de $f$ sur $\Psi(G';\pi)$ se 
r\'eduit \`a l'alg\'ebricit\'e de $f$ sur les vari\'et\'es
$\Psi(G';\tau)$ o\`u $\tau$ est une \ett, qui est \'evidente par les
conditions pos\'ees \`a la construction de $f$. 

Donc $f$ est une fonction trace par application du th\'eor\`eme de
Paley-Wiener ([BDK]). Soit $f'$ une fonction sur $G'$ qui correspond
\`a $f$ par ce 
 th\'eor\`eme. 
\begin{lemme}\label{nulite} L'int\'egrale orbitale de $f'$
 est nulle sur les \'el\'ements semisimples r\'eguliers non
 elliptiques de $G'$.
\end{lemme}

{\bf D\'emonstration.} On a vu que, pour tout $\sigma\in
 Grot_{ind}(G')$, on avait 
 $tr\sigma(f')=0$. Le lemme est alors une cons\'equence classique en
 toute caract\'eristique (voir par exemple
 [Ba1], lemme 2.4).\qed
\ \\

La fonction $f'$ annule de
 plus les traces de toutes les repr\'esentations essentiellement de
 carr\'e int\'egrable, mais pas la trace de $\tau_1$. 
 Montrons qu'il y a l\`a une contradiction qui prouve que $s=1$
 (on est toujours dans la d\'emonstration du lemme \ref{s=1}). On
 distingue deux cas, selon la caract\'eristique de $F$~:\\ 
\ \\
{\bf (a) $F$ est de caract\'eristique nulle}\\

\begin{lemme}\label{nulell} L'int\'egrale orbitale de $f'$ s'annule sur
les \'el\'ements elliptiques r\'eguliers de $G'$. 
\end{lemme}

{\bf D\'emonstration.} Soient $Z$ le centre de $G'$ et $\omega$
un caract\`ere unitaire de $Z$. Nous allons utiliser l'espace
$L_0(G'_e;\omega)$ des fonctions $f$ localement constantes sur l'ensemble
des \'el\'ements elliptiques r\'eguliers de $G'$, invariantes sous
l'action de
conjugaison par des \'el\'ements de $G'$, et v\'erifiant $f(z
g)=\omega (z) f(g)$ pour tout $z\in Z$ et tout $g$ elliptique
r\'egulier dans $G'$. On consid\`ere le sous-espace $L^2(G'_e;\omega)$
de $L_0(G'_e;\omega)$ form\'e des fonctions $f$ telles que
$$\sum_{T\in {\mathcal T}_e} |W_T|^{-1} \int_{T^{reg}/Z}
D(\bar{t})|f(\bar{t})|^2 d\overline{t}$$ 
 converge, o\`u ${\mathcal
T}_e$ est un ensemble de repr\'esentants des classes de conjugaison de
tores elliptiques maximaux de $G'$, $|W_T|$ est le
cardinal du groupe de Weyl de $T$, $T^{reg}$ est l'ensemble des
\'el\'ements r\'eguliers de $T$,  $d\overline{t}$ est choisie de
fa\c{c}on \`a ce que le volume de $T^{reg}/Z$ soit 1, et $D(\bar t)$
est la valeur absolue
normalis\'ee du d\'eterminant de l'op\'erateur  $Ad(t^{-1})-Id$
agissant sur $Lie(G')/Lie(T)$. On d\'efinit un produit scalaire dans
$L^2(G'_e;\omega)$ en 
posant~: 
             $$<f_1;f_2>_e=\sum_{T\in {\mathcal T}_e}
|W_T|^{-1}\int_{T^{reg}/Z} 
D(\bar{t})f_1(\bar{t})\overline{f_2(\bar{t})} d\overline{t},$$ 
qui munit $L^2(G'_e;\omega)$ une structure d'espace
pr\'ehilbertien. 

On sait que pour toute repr\'esentation $\pi$ de $G'$ de carr\'e
int\'egrable et de caract\`ere central $\omega$, la restriction du
caract\`ere $\chi_{\pi}$ de $\pi$ \`a $G'_e$ se trouve dans
$L^2(G'_e;\omega)$ et 
les \'el\'ements de $L^2(G'_e;\omega)$ ainsi obtenus forment une 
famille orthonormale pour $<\ ; >_e$ ([Cl2]). Une
cons\'equence de la correspondance de Jacquet-Langlands avec une
alg\`ebre \`a division et que ce
syst\`eme est complet ([Ba2], cor. 5.13, par exemple).
\`A partir de $f'$ on d\'efinit $f'_{\omega}$ en posant, pour tout
$g\in G'$,
$$f'_{\omega}(g)=\int_Z\omega(z)f'(zg)dz.$$
Nous avons la relation entre les int\'egrales orbitales
$$\Phi(f'_{\omega};g)=\int_Z\omega(z)\Phi(f';zg)dz,$$
pour tout $g\in G'$. Le lemme \ref{nulite} entraine donc que
l'int\'egrale orbitale de 
$f'_{\omega}$ est nulle sur les \'el\'ements semisimples r\'eguliers
qui ne sont pas elliptiques, i.e. qui n'appartiennent pas \`a un tore
elliptique maximal. Dans ce cas, la formule d'int\'egration
de Weyl donne, pour
toute \care\ $\pi$ de $G'$ de caract\`ere central $\omega$, 
$$tr\pi(f)=\sum_{T\in {\mathcal T}_e} |W(T)|^{-1}\int_{T^{reg}/Z}
D(\bar{t})\chi_{\pi}(\bar{t})\Phi(f'_{\omega};\bar{t})d\bar{t}.$$
Or, on sait que $tr\pi(f'_{\omega})=tr\pi(f')=0$.

\begin{lemme} La restriction \`a l'ensemble des \'el\'ements
elliptiques r\'eguliers de $G'$ de l'int\'egrale orbitale de
$f'_{\omega}$ se trouve dans l'espace $L^2(G'_e;\omega^{-1})$.
\end{lemme}

{\bf D\'emonstration.} D'apr\`es [H-CvD], th.14, chap.8, pour tout $T$
dans ${\mathcal T}_e$, le produit de l'int\'egrale orbitale de
$f'_{\omega}$ et de la fonction $D^{\frac{1}{2}}$ est born\'e sur
$T^{reg}$. Le module du carr\'e de cette fonction est ainsi born\'e et
donc int\'egrable sur $T^{reg}/Z$ qui est de mesure 1. Comme
${\mathcal T}_e$ est fini, le r\'esultat s'ensuit. \qed
\ \\

D'apr\`es ce qui pr\'ec\`ede, la conjugu\'ee de $\Phi(f'_{\omega};
\, ^.)$ est un \'el\'ement de $L^2(G'_e;\omega)$ orthogonal \`a
$\chi_{\pi}$. 
C'est vrai pour toute repr\'esentation $\pi$
de carr\'e int\'egrable et de caract\`ere central $\omega$. On en
d\'eduit que $\Phi(f'_{\omega};\, ^.)$ est identiquement nulle sur
l'ensemble des \'el\'ements elliptiques r\'eguliers.
 
Soit $g$ un \'el\'ement elliptique r\'egulier de $G'$. Nous avons
trouv\'e que pour tout caract\`ere unitaire $\omega$ de
$Z$ on a
$$\int_Z\omega(z)\Phi(f';zg)dz=0.$$
Par [Bo], chap. 2, th. 4.4, on en d\'eduit que $\Phi(f';g)=0$. Le
lemme est montr\'e.\qed

Les lemmes \ref{nulite} et \ref{nulell} impliquent que 
l'int\'egrale orbitale de $f'$ est nulle sur l'ensemble des
\'el\'ements semisimples r\'eguliers de $G'$.  
Mais alors la  trace de toute
repr\'esentation de $G'$ est nulle sur $f'$ (nous sommes en
caract\'eristique nulle). Cela  contredit
$tr(\tau_1(f'))\neq 0$. En conclusion $s=1$ et on a
$ind_{P_0}^{G'}\pi_0=a\tau$ o\`u $a$ est un entier strictement positif
et $\tau$ est une \temp\ de $G'$.\qed 
\ \\
\ \\
{\bf (b) $F$ est de caract\'eristique positive}\\

Nous avons utilis\'e deux r\'esultats connus pour l'instant
uniquement en caract\'eristique nulle dans la d\'emonstration plus haut.
Dans le cas de caract\'eristique non nulle, nous allons utiliser les
corps proches de Kazhdan pour conclure. Soit $E$ un
corps local non archim\'edien de caract\'eristique nulle. Notons $O_F$
(resp. $O_E$) 
l'anneau des entiers de $F$ (resp. $E$), et $I_F$ (resp. $I_E$)
l'id\'eal maximal de $O_F$ (resp. $O_E$). Nous disons que $E$ et $F$
sont $m$-proches s'il existe un isomorphisme d'anneaux $\lll$ de $O_F/I_F^m$
sur  $O_E/I_E^m$ (pour tout $m$ on peut trouver un tel corps
$E$). Quand on dira par la suite que $F$ et $E$ sont $m$ proches, on
consid\`erera tacitement qu'un isomorphisme $\lll$ est fix\'e une fois
pour toutes. Soit $E$ un corps $m$ proche de $F$. Rebaptisons $G'_F$
notre groupe pour rappeler qu'il est 
d\'efini sur $F$. Nous avons $G'_F=GL_r(D_F)$ o\`u $D_F$ est une
alg\`ebre \`a division centrale sur $F$, de dimension finie $d^2$. On
note $G'_E$ le groupe 
$GL_r(D_E)$ o\`u $D_E$ est une alg\`ebre \`a division centrale sur $E$
qui a le m\^eme invariant que $D_F$. 

Notons $O_{D_F}$
l'anneau des entiers de $D_F$, et $I_{D_F}$
l'id\'eal maximal de $O_{D_F}$. On pose $K_{D_F}^0=GL_r(D_F)$, et, pour
tout entier strictement positif $l$,
$K_{D_F}^l=1+M_r(I_{D_F}^{dl})$. On fixe sur $G'_F$ une mesure de Haar
telle que le volume de $K_{D_F}^0$ soit 1. Soit $f$ une fonction localement
constante \`a support compact sur $G'_F$. Nous d\'efinissons le {\it niveau}
de $f$ comme \'etant le plus petit entier $l$ tel que $f$ soit
bi-invariante par $K_{D_F}^l$. Notons $H_F^l$ l'alg\`ebre de Hecke des
fonctions localement constantes \`a support compact de niveau
inf\'erieur ou \'egal \`a $l$ sur $G'_F$. Si $\pi$ est une
repr\'esentation lisse     
irr\'eductible de $G'_F$ on appelle niveau de $\pi$ le plus petit
entier $l$ tel que $\pi$ ait un vecteur non nul fixe sous
$K_{D_F}^l$. Adoptons pour $G'_E$ des notations et d\'efinitions
analogues \`a celles 
fix\'ees dans ce paragraphe pour $G'_F$.

Dans [Ba2] nous montrons que, quel que soit l'entier
positif $l$, il existe un entier $m$ tel que, pour tout corps $E$
$m$-proche de $F$, $\ll$ induise un isomorphisme  d'alg\`ebres $\zzd$
de $H_F^l$ sur $H_E^l$. D'o\`u une bijection entre
l'ensemble de classes d'\'equivalence des repr\'esentations lisses
irr\'eductibles de $G'_F$ de niveau inf\'erieur ou \'egal \`a $l$ et  
l'ensemble de classes d'\'equivalence des repr\'esentations lisses
irr\'eductibles de $G'_E$ de niveau inf\'erieur ou \'egal \`a $l$. On
utilise pour  cet isomorphisme la m\^eme notation $\zzd$. Nous avons 
$$tr \zzd (\pi)(\zzd(f))=tr\pi(f)$$
pour toute repr\'esentation $\pi$ de $G'_F$ de niveau inf\'erieur ou
\'egal \`a $l$ et toute $f\in H_F^l$.

Avec les conventions faites en d\'ebut de
section, un \sgls\ $L$ de $G'_F$ ou $G'_E$ est form\'e des
matrices diagonales par blocs de taille donn\'ee et on peut associer
\`a $L$ de fa\c{c}on biunivoque une suite ordonn\'ee d'entiers
strictement positifs $r_1,r_2...r_k$ telle que $\sum_{i=1}^kr_i=r$
o\`u les $r_i$ repr\'esentent les tailles de ces blocs. \`A un \sgls\
de $G'_F$ correspond donc un unique \sgls\ de $G'_E$. C'est pareil
pour les sous-groupes paraboliques standard. Si $L_F$ est un \sgls\ de
$G'_F$ on note $L_E$ le \sgls\ de $G_E$ qui lui correspond, et si
$P_F$ est un sous-groupe parabolique standard de $G'_F$ on note $P_E$
le sous-groupe parabolique standard de $G'_E$ qui lui correspond. Sur
chaque bloc d'un \sgls\ de $G'_F$ et $G'_E$ nous adoptons les m\^emes
notations et conventions que plus haut. \\
\ \\

Soit $P$ est un sous-groupe parabolique standard de $G'_F$ et $P=LU$ sa
d\'ecomposition de Levi ($L$ \sgls\ de $G'_F$ et $U$ le radical
unipotent de $P$). On munit $P$ et $U$
de mesures de Haar invariantes \`a gauche telles que les volumes de
$P\cap K_{D_F}^0$ et $U\cap K_{D_F}^0$ soient \'egaux \`a 1.
On note
 $\delta_P$ le caract\`ere modulaire sur $P$.  

Pour toute fonction localement constante \`a
support compact sur $G'_F$, on d\'efinit une fonction localement
constante \`a support constant sur $L$, $f^P$, par la formule~:
$$f^P(l)=\delta_P^{\frac{1}{2}}(l)\int_{K_{D_F}^0}\int_Uf(k^{-1}luk)dkdu$$
pour tout $l\in L$.

Les m\^emes d\'efinitions
s'appliquent \`a $E$.

\begin{prop}\label{parab} Soient $P_F$ un sous-groupe parabolique
standard de $G'_F$ et $P_F=L_FU_F$ la d\'ecomposition de Levi standard
de $P_F$. Alors il existe un $m$ tel que, si $F$ et $L$ sont
$m$-proches, alors $\zzd(f')$ et $\zzd(f^{'P_F}$) sont bien d\'efinies
et on a  
$$(\zzd(f'))^{P_E}=\zzd(f^{'P_F}).$$ 
\end{prop}

{\bf D\'emonstration.} L'analogue de cette proposition  dans le cas
particulier $G'_F=GL_r(F)$ est montr\'e \`a la page 1053 de [Le]. La
d\'emonstration est exactement la m\^eme dans notre cas.\qed
\ \\

Nous allons relever notre situation en caract\'eristique nulle et
utiliser (a).
Soit ${\mathcal P}_F$ l'ensemble de tous les sous-groupes paraboliques
standard, propres ou pas, de $G'_F$ et soit $m$ un entier suffisamment
grand pour que, pour tout corps local $E$ qui est $m$-proche de $F$,
la proposition \ref{parab} plus haut soit v\'erifi\'ee  pour tout
$P_F\in {\mathcal P}_F$ qui a une d\'ecomposition de Levi standard
$P_F=L_FU_F$ (c'est possible  puisque ${\mathcal P}_F$ est un ensemble
fini). On a alors~:\\ 

1) si $P_E$ est un  sous-groupe parabolique standard propre de $G'_E$,
   si $P_E=L_EU_E$ est une d\'ecomposition de Levi standard de $P_E$,
   si $\pi$ est une repr\'esentation lisse irr\'eductible de $L_E$,
   alors~:  

\ \ \ \ - ou bien le niveau de $\pi$ est sup\'erieur strictement \`a
$l$ et alors $$tr\pi(\zzd(f^{'P_F}))=0,$$ 

\ \ \ \ - ou bien le niveau de $\pi$ est inf\'erieur ou \'egal  \`a
$l$ et alors, en supposant que $\sigma$ est la repr\'esentation lisse
irr\'eductible de $L_F$ qui v\'erifie $\zzd(\sigma)=\pi$ on peut
\'ecrire~:  
$$tr\pi(\zzd(f^{'P_F}))=tr\sigma(f^{'P_F})=tr (ind_{P_F}^{G'_F}\sigma)
(f')=0$$ car  $f'$ annule la trace de toute repr\'esentation
de $Grot_{ind}(G'_F)$. 

Dans les deux cas on obtient $tr\pi(\zzd(f^{'P_F}))=0$ ce qui
implique, compte tenu de  la proposition \ref{parab} plus haut, que  
$$tr(ind_{P_E}^{G'_E}\pi)(\zzd(f'))=0.$$ 
 Finalement, on a trouv\'e que $\zzd(f')$ est une fonction qui annule
la trace de toute repr\'esentation dans $Grot_{ind}(G'_E)$.\\ 

2) si $\pi$ est une \eci\ de $G_E$, alors ou bien son niveau est
   sup\'erieur strictement \`a 
   $l$ et donc $tr\pi(\zzd(f'))=0$ ou bien son niveau est inf\'erieur
   ou \`egal \`a $l$ et alors, en supposant que $\sigma$ est la
   repr\'esentation lisse irr\'eductible de $G'_F$ qui v\'erifie
   $\zzd(\sigma)=\pi$ on a que $\sigma$ est une \eci\
   (th. 2.16.b dans [Ba2]) et par cons\'equent  
 $$tr\pi(\zzd(f'))=tr\sigma(f')=0,$$ 
car $f'$ annule la trace de toute \eci\ de $G'_F$. Finalement, on a
trouv\'e que $\zzd(f')$ annule la trace de toute \eci\ de $G'_E$.\\ 

Or, $E$ est de caract\'eristique nulle.
Par les points 1) et 2) ci-dessus et par le raisonnement d\'ej\`a fait 
en caract\'eristique nulle, on en d\'eduit que $\zzd(f')$ annule la
trace de toutes les repr\'esentations de $G'_E$. D'autre part,  comme
$tr\tau_1(f')\neq 0$, le niveau de $\tau_1$ est inf\'erieur ou \'egal
\`a $l$ et donc $\zzd(\tau_1)$ est bien d\'efini et  
$$tr\zzd(\tau_1)(\zzd(f'))=tr\tau_1(f')\neq 0.$$
 Contradiction. On trouve donc, comme en caract\'eristique nulle, que
$s=1$ et on a $ind_{P_0}^{G'}\pi_0=a\tau$ o\`u $a$ est un entier
strictement positif et $\tau$ est une \temp\ de $G'$.\qed 

Le fait que $a=1$ se montre exactement de la m\^eme fa\c{c}on que les
\'etapes (2) et (3) de la preuve de la proposition 27 de [FK], page
98. Le th. \ref{irred} est d\'emontr\'e.

\section{ (Une) correspondance de Jacquet-Langlands g\'en\'eralis\'ee} 

\subsection{ Correspondance entre les groupes de Grothendieck}\label{grot}

Dans cette sous-section $F$ est un corps local non archim\'edien, $D$
est une alg\`ebre \`a division centrale sur $F$ de dimension $d^2$, $r$ est
un entier strictement positif, $n=rd$, $G=GL_n(F)$ et
$G'=GL_r(D)$.
Rappelons que, si $g$ est un \'el\'ement de $G$ ou de $G'$ on dit que
$g$ est {\it semisimple r\'egulier} si son polyn\^ome
caract\'eristique est s\'eparable (valeurs propres distinctes dans une
cl\^oture alg\'ebrique de $F$). Si $g\in G$ et $g'\in G'$  on dit que
$g$ et $g'$ {\it se correspondent} s'ils sont semisimples r\'eguliers et
ont le m\^eme polyn\^ome caract\'eristique. Rappelons maintenant
l'\'enonc\'e de la
correspondance de Jacquet-Langlands~: 
{\it il existe une unique
correspondance bijective $\bf C$ de l'ensemble des classes d'\'equivalence de
\eci\ de $G$ sur l'ensemble des classes d'\'equivalence de
\eci\ de $G'$ qui v\'erifie, pour toute \eci\ $\pi$ de $G$~:
$$\chi_{\pi}(g)=(-1)^{n-r}\chi_{{\bf C}(\pi)}(g')$$
pour tout $g$ et $g'$ qui se correspondent.}
Sur $G$ comme sur $G'$ fonctionnent les conventions
faites \`a la section pr\'ec\'edente pour une forme int\'erieure de
$G$ quelconque. Avec ces conventions, \`a un \sgls\ de $G'$
correspond un unique \sgls\ de $G$ et \`a un sous-groupe parabolique
standard de $G'$ correspond un unique sous-groupe parabolique standard
de $G$. Si $L'$ est un \sgls\ de $G'$ et $L$ est le \sgls\ de $G$ qui
lui correspond, $L'$ est un produit de groupes du type $GL_k(D)$ et
$L$ est le produit des groupes du type $GL_{dk}(F)$ qui leur
correspondent. Si $L$ est de ce type, on dit que $l$ {\it se
transf\`ere}. Ces d\'efinitions s'\'etendent facilement aux
sous-groupes paraboliques standard.
La 
correspondance de Jacquet-Langlands induit une correspondance
 entre les \ecis\ de $L$
et les \ecis\ de 
$L'$. Nous noterons cette correspondance aussi {\bf C}, quand il n'y
aura pas d'ambigu\"\i t\'e.\\

\begin{prop}\label{base} Soit ${\mathcal B}_G$ le sous-ensemble de
$Grot(G)$ form\'e des induites de \ecis\ de tous les sous-groupes
paraboliques standard de $G$. Alors ${\mathcal B}_G$ n'est autre que
la base de Langlands de $Grot(G)$. Soit ${\mathcal B}_{G'}$ le
sous-ensemble de $Grot(G')$ form\'e des induites de \ecis\ de tous les
sous-groupes 
paraboliques standard de $G'$. Alors ${\mathcal B}_{G'}$ n'est autre
que la base de Langlands de $Grot(G')$. 
\end{prop}

{\bf D\'emonstration.} Le fait que ${\mathcal B}_G$ soit une base de
$Grot(G)$ appara\^\i t dans [Ze]  sans que le fait que ce soit exactement la
base de Langlands y soit dit explicitement. Ce r\'esultat est
cons\'equence directe du fait que les induites des \cares\ sont 
irr\'eductibles et la d\'emonstration est la m\^eme sur les deux
groupes. Voir [Ta], prop. 2.1 pour le cas g\'en\'eral de $G'$.\qed

Soit $\jl$ l'application d\'efinie sur ${\mathcal B}_{G'}$ \`a valeurs
dans $Grot(G)$ de la fa\c{c}on suivante~:  
soient $P'$ un sous-groupe parabolique standard (propre ou pas) de
$G'$, $L'$ le sous-groupe de Levi standard de $P'$ et $\sigma$ une
\eci\ de $L'$~; soient $P$ le sous-groupe  
parabolique standard de $G$ qui correspond \`a $P'$, $L$ le
sous-groupe de Levi standard de $P$ ($L$ correspond \`a $L'$) et
$\sigma$ la \eci\ de $L$ qui v\'erifie ${\bf C}(\sigma)=\sigma'$~; on
pose  
$$\jl(ind_{P'}^{G'}\sigma')=ind_{P}^{G}\sigma.$$
 Gr\^ace \`a la proposition pr\'ec\'edente, $\jl$ s'\'etend de
fa\c{c}on unique en un morphisme injectif de ${\mathbb{Z}}$-modules de
$Grot(G')$ dans $Grot(G)$. 

\begin{theo}\label{jlr} a) Le morphisme injectif de ${\mathbb{Z}}$-modules
$$\jl~: Grot(G')\to Grot(G),$$ 
commute \`a la tensorisation avec des caract\`eres, et on a 
$$\chi_{\pi'}(g')=(-1)^{n-r}\chi_{\jl (\pi')}(g)$$ 
pour tout $g\in G$ et $g'\in G'$ qui se correspondent.

b) Le morphisme $\jl$ r\'ealise une bijection 

- entre l'ensemble des classes de \cares\ de $G'$ et l'ensemble de
  classes de \cares\ de $G$, ainsi qu'entre l'ensemble des classes de 
  \ecis\ de $G'$ et l'ensemble des classes de \ecis\ de $G$
  (l'application inverse n'est 
  autre que ${\bf C}$), 

- entre l'ensemble des classes de 
\temps\ de $G'$ et l'ensemble des classes de  \temps\ de $G$ qui sont
des 
  induites de \cares\ de sous-groupes paraboliques de $G$ qui se
  transf\`erent, 

- entre l'ensemble des classes de repr\'esentations cuspidales de $G'$
  et l'ensemble de classes de \ecis\  de $G$ 
  dont la restriction de Jacquet \`a tout sous-groupe parabolique qui
  se transf\`ere est nulle. 

c) Soit $\ss$ le sous-module de $Grot(G)$ engendr\'e par l'ensemble
des repr\'esentations induites de \ecis\ de tous les sous-groupes
paraboliques de $G$ qui ne se transf\`erent pas. Alors $\ss$ est un
suppl\'ementaire de l'image $Im(\jl)$ de $\jl$ dans $Grot(G)$ et $\jl$
induit un isomorphisme de groupes
$$\overline{\jl} : Grot(G')\simeq Grot(G)/\ss.$$
\end{theo}

{\bf D\'emonstration.} a) Soient $g\in G$ et $g'\in G'$ qui se
correspondent. Pour v\'erifier qu'on a

$$\chi_{\pi'}(g')=(-1)^{n-r}\chi_{\jl (\pi')}(g)\ \ \ \ \ \ \ \
\forall\ g\lra g'$$ 
\ \\
pour tout $\pi'\in Grot(G')$ il suffit de le v\'erifier pour tout
$\pi'\in {\mathcal B}_{G'}$. Si $\pi'$ est une \eci\ on a
$\jl(\pi')={\bf C}^{-1}(\pi')$ et c'est \'evident. Si
$\pi'$ est l'induite d'une \eci\ \`a partir d'un sous-groupe
parabolique propre, cela r\'esulte de la
d\'efinition de $\jl(\pi')$ et de la formule du caract\`ere d'une
repr\'esentation induite ([vD] ou [Cl1], prop. 3). 

Pour montrer que $\jl$ commute \`a la tensorisation avec des
caract\`eres il suffit de le montrer sur la base ${\mathcal
B}_{G'}$. Soit $\chi$ un caract\`ere de   $G'$. Si
$\pi'=ind_{P'}^{G'}\sigma'$, o\`u $P'$ est un sous-groupe parabolique
standard de $G'$ et $\sigma'$ est une \eci\ du sous-groupe de Levi
standard $L'$ de  $P'$, alors on a  
$$\chi\otimes \pi'=ind_{P'}^{G'}((res_{P'}^{G'}\chi)\otimes\sigma').$$
Or, $(res_{P'}^{G'}\chi)\otimes\sigma'$ est une \eci\ de $L'$. Soient
$P$ le sous-groupe parabolique standard de $G$ correspondant \`a $P'$,
$L$ son sous-groupe de Levi standard et $\sigma$ la \eci\ de $L$ qui
correspond \`a $\sigma'$ par la correspondance entre $L$ et $L'$. On a
alors par d\'efinition de $\jl$~: 
$$\jl(ind_{P'}^{G'}((res_{P'}^{G'}\chi)\otimes\sigma'))=
ind_{P}^{G}((res_{P}^{G}\chi)\otimes\sigma).$$   
On en d\'eduit que 
$$\jl(\chi\otimes\pi')=\jl(ind_{P'}^{G'}((res_{P'}^{G'}\chi)\otimes\sigma'))=
ind_{P}^{G}((res_{P}^{G}\chi)\otimes\sigma)= 
$$
\\  
$$\chi\otimes
ind_P^G\sigma=\chi\otimes\jl(ind_{P'}^{G'}\sigma')=\chi\otimes
\jl(\pi').$$ 
\ \\

b) L'assertion sur les \ecis\ est la correspondance de
Jacquet-Langlands classique. L'assertion sur les \cares\ en r\'esulte
puisque {\bf C} conserve le caract\`ere central et une \eci\ est de carr\'e
int\'egrable si et seulement si son caract\`ere central est
unitaire. Les \temps\ de $G'$ sont exactement les induites de \cares\
des sous-groupes paraboliques (propres ou pas) et l'assertion sur les
\temps\ est imm\'ediate. La partie concernant les \cusps\ est le
th\'eor\`eme B.2.b 
de [DKV] (l'argument marche aussi en caract\'eristique non nulle, voir
[Ba2], prop. A, page 53).\\ 

c) Le fait que $\ss$ est un sous-module suppl\'ementaire de $Im(\jl)$
dans $Grot(G)$ est une cons\'equence imm\'ediate du fait que
$\jl({\mathcal B}_{G'})\subset {\mathcal B}_G$ et 
${\mathcal B}_G\ \backslash\  \jl({\mathcal B}_{G'})$ 
est par d\'efinition une base
de $\ss$. Le fait que le morphisme de ${\mathbb Z}$-modules
$Grot(G')\to Grot(G)/ \ss$ obtenu par composition de la projection
$Grot(G)\to Grot(G)/ \ss$ avec $\jl$ est un isomorphisme s'ensuit,
parce que $\jl$ est injectif et r\'ealise donc une bijection de
$Grot(G')$ sur son image.\qed
\ \\
\ \\
{\bf Commentaires.}

1) Le morphisme  $\jl$ n'envoie pas toute repr\'esentation cuspidale
   de $G'$ sur une repr\'esentation cuspidale de $G$. C'est le cas si
   et seulement si l'ensemble des nombres premiers qui divisent $n$
   co\"\i ncide avec l'ensemble des nombres premiers qui divisent $r$.\\ 

2) Le morphisme $\jl$  n'envoie pas toute repr\'esentation
   irr\'eductible de $G'$ sur une repr\'esentation irr\'eductible de
   $G$, sauf si $d=1$ ou $r=1$. Donons un contre-exemple qui se
   g\'en\'eralise facilement. Soient $n=4$, $d=r=2$, $L'$ le tore
   diagonal de $G'=GL_2(D)$ et $L$ le \sgls\ de $G=GL_4(F)$ qui lui
   correspond. Soit $\rho$ un caract\`ere de $F$. En utilisant les
   notations de [Ze], consid\'erons la \eci\ $\sigma$ de $L$~: 
$$\sigma=<[\rho;\nu\rho]>^t\otimes <[\nu\rho;\nu^2\rho]>^t$$
($\nu$ est le caract\`ere valeur absolue du d\'eterminant,
$<[\rho;\nu\rho]>^t$ est la r.e.c.i. associ\'ee au segment
$[\rho;\nu\rho]$ etc.).
Le th\'eor\`eme 9.7 de [Ze] implique que la repr\'esentation induite de
$\sigma$ \`a $G$ est r\'eductible, puisque les segments
$[\rho;\nu\rho]$ et $[\nu\rho;\nu^2\rho]$ sont li\'es. Soit $\sigma'$
la \eci\ de $L'$ qui correspond \`a $\sigma$. Alors le lemme 2.5 de
[Ta] implique que la repr\'esentation induite de $\sigma'$ \`a $G'$
est irr\'eductible parce que $d=2$ ne divise pas la longueur du segment
$[\rho;\nu\rho]\cap [\nu\rho;\nu^2\rho]=[\nu\rho;\nu\rho]$.\\ 

On peut donc dire qu'en g\'en\'eral, \`a une repr\'esentation
irr\'eductible $\pi$ de $G'$ correspond une combinaison lin\'eaire \`a
coefficients dans ${\mathbb Z}$ de repr\'esentations irr\'eductibles
de $G$.\\ 

3) Le th\'eor\`eme \ref{jlr} implique que, si $D$ et $D'$ sont deux
   alg\`ebres \`a division de m\^eme dimension sur $F$, si $r$ est un
   entier strictement positif, il y a un {\it unique} isomorphisme  
$$f~: Grot(GL_r(D))\simeq Grot(GL_r(D'))$$
tel que pour tout $\pi \in Grot(D)$ on ait~:
$$\chi_{\pi}(g)=\chi_{f(\pi)}(g')$$ 
pour tout $g\in GL_r(D)$ et tout $g'\in GL_r(D')$ semisimples
r\'eguliers qui ont le
m\^eme polyn\^ome caract\'eristique. Une question int\'eressante
serait de savoir si cette fois, la correspondance $f$ envoie une
repr\'esentation irr\'eductible sur une repr\'esentation
irr\'eductible. Cela n'a pas l'air \'evident, m\^eme si le
th\'eor\`eme \ref{jlh} de la sous-section suivante implique que la
correspondance $f$ commute avec l'induction et la restriction, ainsi
qu'avec l'involution d\'efinie par A.-M. Aubert ([Au]).\\ 

4) En d\'efinissant le morphisme $\jl$ ainsi on a fait le choix
   (naturel d'un certain point de vue) de privil\'egier les \cares. En
   effet, les propri\'et\'es impos\'ees au point a) du th\'eor\`eme ne
   caract\'erisent pas le morphisme $\jl$, et m\^eme les
   propri\'et\'es impos\'ees au points a) et b) ne le caract\'erisent
   pas non plus. En regardant ``dans l'autre sens'' on  obtient le
   r\'esultat suivant plus naturel d'un autre point de vue puisqu'il
   est assorti d'``unicit\'e''~:\\ 

\begin{prop}\label{antiisomgrot} Il existe un unique morphisme de
   $\mathbb Z$-modules  
$${\bf LJ}_r : Grot(G)\to Grot(G')$$
tel que, pour tout $\pi'\in Grot(G')$, l'image r\'eciproque de $\pi'$
par ${\bf LJ}_r$ est l'ensemble de tous les $\pi\in Grot(G)$ qui
v\'erifient~: 
$$\chi_{\pi}(g)=(-1)^{n-r}\chi_{\pi'}(g')$$
pour tout $g\in G$ et $g'\in G'$ qui se correspondent.
Le morphisme ${\bf LJ}_r$ est surjectif. Son noyau est $\ss$.\\ 
\end{prop}

Disons qu'un \'el\'ement de $Grot(G)$ est $G'-nul$ si son caract\`ere est
nul sur tout \'el\'ement semisimple r\'egulier qui a un correspondant
dans $G'$. L'ensemble des \'el\'ements $G'-nuls$ dans $Grot(G)$ n'est autre 
que $\ss$.

\ \\

\subsection{ Correspondance entre les alg\`ebres de Hopf}

Dans cette sous-section, $F$ d\'esigne un corps local non
archim\'edien et $D$ une alg\`ebre \`a division centrale sur $F$, de
dimension $d^2$.  

Soit $n$ un entier strictement positif. Soit $G$ l'un des groupes
$GL_n(F)$ ou $GL_n(D)$. Rappelons que, si $L$ est un \sgls\ de $G$,
alors il existe une 
partition $A_1\coprod A_2\coprod ...\coprod A_k$ de l'ensemble $\{1,2\
...n\}$ o\`u $A_1=\{1,2\ ...n_1\},\ A_2=\{n_1+1,n_1+2\ ...n_1+n_2\}$
etc. telle que $L$ soit l'ensemble des matrices $M=(m_{ij})_{1\leq
i,j\leq n}\in G$ telles que $m_{ij}$ est nul si $(i;j)\notin
\cup_{i=1}^kA_i\times A_i$. On appelle les ensembles ordonn\'es
$A_1,A_2 ...A_k$ les {\it sections} de $L$ et on pose aussi
$|L|=n-k-1$. Maintenant, si $L_1$ et $L_2$ sont deux sous-groupes de
Levi standard de $G$, on note $W(L_1,L_2)$ le sous-groupe de $W(G)$
form\'e des \'el\'ements $w$ qui v\'erifient~: 

- $w(k)<w(l)$ si $k<l$ et $k$ et $l$ sont dans la m\^eme section de $L_1$

- $w^{-1}(k)<w^{-1}(l)$  si $k<l$ et $k$ et $l$ sont dans la m\^eme
  section de $L_2$. 

Si $w\in W(L_1;L_2)$, alors $^wL_1\cap L_2$ et $L_1\cap L_2^w$ sont
des sous-groupes de Levi standard de $G$. Dans ce qui suit, on
n'utilise les foncteurs d'induction et de restriction qu'\`a partir de
sous-groupes paraboliques standard de $G$ qui vont donc \^etre
d\'etermin\'es par leur sous-groupes de Levi standard. On adopte donc
les notations de [Au]~: si $L_1$ et $L_2$ sont deux sous-groupes de
Levi standard de $G$ tels que $L_1\subset L_2$, alors on pose
$\iii_{L_1}^{L_2}=ind_{P_1}^{P_2}$ et $\r_{L_1}^{L_2}=res_{P_1}^{P_2}$
o\`u $P_1$ et $P_2$ sont les sous-groupes paraboliques standard qui
ont pour sous-groupes de Levi standard $L_1$ et $L_2$
respectivement.\\  

Rappelons ici sous forme d'un lemme le th\'eor\`eme 1.2 de [Ze],
r\'esultat qui nous sera utile par la suite~:

\begin{lemme}\label{geom} Soient $L_1$ et $L_2$ deux sous-groupes de
Levi standard de $G$. Soit $\pi$ une repr\'esentation de $L_1$. On a
alors~: 
$$\r_{L_2}^{G}\iii_{L_1}^G\pi=\sum_{w\in W(L_1,L_2)}\iii_{^wL_1\cap
L_2}\bigl( (\r_{L_1\cap L_2^w}\pi)^w\bigr).$$ 

\end{lemme}

Posons 
$$\rr(F)=\bigoplus_{n\in {\mathbb N}}Grot(GL_n(F))$$ 
o\`u, par convention, pour $n=0$ on met
$GL_n(F)={\mathbb Z}=Grot(GL_n(F))$. Si $L$ est un \sgls\ de
$GL_n(F)$ 
\`a $k$ blocs, on consid\'erera $\r_L^G$ comme une application
lin\'eaire de $Grot(GL_n(F))$ \`a valeurs dans $\displaystyle\otimes^k
\rr(F)$.  

 A. Zelevinski a muni l'espace vectoriel $\rr(F)$ d'une multiplication
 $m$ et une comultiplication $c$ qui en font une alg\`ebre de Hopf
 gradu\'ee ([Ze]). Avec les notations et conventions plus haut, la
 multiplication et la comultiplication de Zelevinski se d\'efinissent
 de la fa\c{c}on suivante~: 

- si $\pi_k\in Grot(GL_k(F))$ et $\pi_{k'}\in Grot(GL_{k'}(F))$, alors
  on identifie $GL_k(F)\times GL_{k'}(F)$ avec le sous-groupe de Levi
  standard $L$ de  $GL_{k+k'}(F)$ de blocs diagonaux de taille $k$
  puis $k'$ et on pose  
$$m(\pi_k;\pi_{k'})=\iii_L^G(\pi_k \otimes \pi_{k'})$$ 
et on \'etend par bilin\'earit\'e \`a une application de
$\rr(F)\otimes \rr(F)$ \`a valeurs dans $\rr(F)$. 

- si $\pi\in Grot(GL_n(F))$, on identifie $GL_k(F)\times GL_{k'}(F)$,
  $k+k'=n$, \`a un sous-groupe de Levi standard maximal de $GL_n(F)$,
  $L_{k,k'}$, et on pose  
$$c(\pi)=\sum_{k+k'=n}\r_{L_{k,k'}}^{GL_n(F)}\pi\ ;$$
on \'etend ensuite par lin\'earit\'e \`a une application de $\rr(F)$
\`a valeurs dans $\rr(F)\otimes \rr(F)$.\\ 

Exactement de la m\^eme fa\c{c}on on peut poser 
$$\rr(D)=\bigoplus_{r\in {\mathbb N}}Grot(GL_r(D)).$$
 Si $L$ est un \sgls\ de $GL_r(D)$ \`a $k$ blocs, on regarde $\r_L^G$
comme une application lin\'eaire de $Grot(GL_r(D))$ \`a valeurs dans
$\otimes^k\rr(D)$. On peut munir l'espace vectoriel $\rr(D)$ d'une
multiplication $m'$ et une comultiplication $c'$ d\'efinies comme plus
haut et v\'erifier comme dans [Ze] qu'on obtient ainsi une alg\`ebre de
Hopf gradu\'ee (voir [Ta]). Si $H$ est l'un des groupes \'etudi\'es
ici, notons  $\Pi^2(H)$ l'ensemble des classes d'\'equivalence de
r.e.c.i. de $H$. 

\begin{prop}\label{inele} En tant qu'anneau, $\rr(F)$ est isomorphe
\`a l'anneau de polyn\^omes en une infinit\'e de variables
commutatives ${\mathbb{Z}}[\displaystyle\bigcup_{n\in
{\mathbb{N}}}\Pi^2(GL_n(F))]$ et, toujours en tant qu'anneau, $\rr(D)$
est isomorphe \`a l'anneau de polyn\^omes en une infinit\'e de
variables commutatives ${\mathbb{Z}}[\displaystyle\bigcup_{r\in
{\mathbb{N}}}\Pi^2(GL_r(D))].$ 
\end{prop}

{\bf D\'emonstration.} Pour $\rr(F)$, cette proposition est une
cons\'equence du corrolaire 7.5, du th\'eor\`eme 9.3, et de la
proposition 9.16 de [Ze]. Voir [Ta], prop. 2.1 pour le cas g\'en\'eral
$\rr(D)$. Autant pour $\rr(F)$ que pour $\rr(D)$ c'est
une cons\'equence imm\'ediate de la proposition \ref{base}.\qed 

\begin{theo}\label{jlh}  a) L'ensemble des  morphismes injectifs de
${\mathbb Z}$-modules  
$$\jl:Grot(GL_r(D))\to Grot(GL_{rd}(F))$$
pour tout $r$ induit un morphisme injectif d'\ {\rm anneaux} 
$$\jla : \rr(D)\to \rr(F).$$

b) L'image de $\jla$ est le sous-anneau de $\rr(F)$ engendr\'e par les
\ecis\ de tous les $GL_n(F)$ tels que $d$ divise $n$. L'id\'eal ${\bf
I}_{F,D}$ de l'anneau $\rr(F)$ engendr\'e par les \ecis\ de tous les
$GL_n(F)$ tels que $d$ ne divise pas $n$ est un ${\mathbb Z}$-module
suppl\'ementaire de l'image de $\jla$ dans $\rr(F)$.  

c) La comultiplication  de $\rr(F)$ induit une
op\'eration bien d\'efinie sur l'anneau $\rr(F)/\ii$. L'anneau
$\rr(F)/\ii$ h\'erite ainsi d'une structure d'alg\`ebre de
Hopf. L'application $\jla$ induit  un isomorphisme d'\ {\rm 
alg\`ebres de Hopf}  
$$\jlh : \rr(D)\simeq \rr(F)/{\bf I}_{F,D}.$$ 
\end{theo}

{\bf D\'emonstration.} a) Pour tout entier positif $r$, $\jl$ induit
une bijection de 
$\Pi^2(GL_{r}(D))$ sur $\Pi^2(GL_{rd}(F))$.  La proposition
\ref{inele} implique que cette restriction de $\jl$ aux \ecis\ pour
tout $r$ induit un unique morphisme d'anneaux $\jla$ de $\rr(D)$ dans
$\rr(F)$, qui de plus est injectif. On doit prouver que pour tout $r$
l'application $\jla$ co\"\i ncide avec $\jl$ sur $Grot(GL_r(D))$ tout
entier. Par la proposition \ref{base} et par lin\'earit\'e des deux
applications, il suffit de v\'erifier que pour deux \ecis\, $\pi'_1$
de $GL_{r_1}(D)$ et $\pi'_2$ de $GL_{r_2}(D)$, on a  
$$m(\jla_{r_1}(\pi'_1);\jla_{r_2}(\pi'_2))=
\jla_{r_1+r_2}(m'(\pi'_1;\pi'_2)).$$
 Mais cela fait partie de la d\'efinition m\^eme de
$\jla_{r_1+r_2}$. \\ 

b) Partant de la d\'efinition de $\jla$, le fait que l'image de $\jla$
est le sous-anneau de $\rr(F)$ engendr\'e par les \ecis\ de tous les
$GL_n(F)$ tels que $d$ divise $n$ est tautologique. En \'ecrivant
maintenant  
$$\bigcup_{n\in {\mathbb{N}}}\Pi^2(GL_n(F))=\big(
\bigcup_{d|n}\Pi^2(GL_n(F))\big)\bigcup \big( 
\bigcup_{d \not{\, \mid}n}\Pi^2(GL_n(F))\big),$$ 
on a que tout polyn\^ome en les variables $\bigcup_{n\in
{\mathbb{N}}}\Pi^2(GL_n(F))$ s'\'ecrit de fa\c{c}on unique comme la
somme d'un polyn\^ome en les variables $\bigcup_{d|n}\Pi^2(GL_n(F))$
et un polyn\^ome dont chaque mon\^ome contient au moins une variable
se trouvant dans l'ensemble $\bigcup_{d \not{\, \mid} n}\Pi^2(GL_n(F))$. On en
d\'eduit l'\'egalit\'e de ${\mathbb{Z}}$-modules~: 
$$\rr(F)=Im(\jla)\oplus \ii.$$

c) Par le point b) le morphisme injectif d'anneaux $\jla:\rr(D)\to
\rr(F)$ induit un isomorphisme d'{\it anneaux} $\jlh:\rr(G')\to
\rr(F)/\ii$. Il suffit de montrer que la comultiplication $c$
 de $\rr(F)$ ``passe au quotient'' et que $\jlh$
commute \`a  la comultiplication.  
Reprenons les notations de la sous-section pr\'ec\'edente~: $r\in
{\mathbb N}^*$, $G'=GL_r(D)$ et $G=GL_{rd}(F)$. On va \'etudier
l'effet du foncteur restriction sur les bases ${\mathcal B}_G$ et
${\mathcal B}_{G'}$. On commence par les \ecis.  

 Dans ce qui suit on suppose que le lecteur est un peu familiaris\'e
 avec les travaux de Bernstein et Zelevinski pr\'esent\'es dans [Ze].\\ 

\begin{lemme}\label{0.2.7} a) Soit $\sigma$ une \eci\ de $G$ qui
correspond \`a un segment $[\rho;\nu^{k-1}\rho]$ de Zelevinski. Soit
$L$ un sous-groupe de Levi standard de $G$. Alors $\r_L^G\sigma$ est
la repr\'esentation nulle si l'entier $n/k$ ne divise pas la taille de
tous les blocs de $L$ et $\r_L^G\sigma$ est une \eci\ de $L$ si
l'entier $n/k$ divise la taille de tous les blocs de $L$. 
 
b) Soient $\sigma'$ une \eci\ de $G'$ et $L'$ un \sgls\ de
$G'$. Posons $\sigma=\jl(\sigma')$. Supposons que la \eci\ $\sigma$ de
$G$ correspond \`a un segment $[\rho;\nu^{k-1}\rho]$. Soit $L$ le
\sgls\ de $G$ qui correspond \`a $L$. Alors $\r_{L'}^{G'}\sigma'$ est
la repr\'esentation nulle si l'entier $n/k$ ne divise pas la taille de
tous les blocs de $L$ et $\r_{L'}^{G'}\sigma'$\ est une \eci\ de $L'$
si l'entier $n/k$ divise la taille de tous les blocs de $L$. 
$\r_{L'}^{G'}\sigma'$ et $\r_L^G\sigma$ se correspondent par {\bf C}. 
\end{lemme}

{\bf D\'emonstration.} Le point a) est la proposition 9.5 de [Ze]. Le
point b) est alors la cons\'equence directe de  [DKV], th. B.2.b.\qed

Appliquons maintenant le lemme \ref{geom} \`a un \'el\'ement $\pi$ de
${\mathcal B}_G$. On a $\pi=\iii_{L_1}^G\sigma$ o\`u $\sigma$ est une
\eci\ du \sgls\ $L_1$ de $G$. Soit $L_2$ un autre \sgls\ de $G$. On a
alors  
$$\r_{L_2}^G\pi=\r_{L_2}^G\iii_{L_1}^G\sigma=\sum_{w\in
W(L_1,L_2)}\iii_{^wL_1\cap L_2}\big( (\r_{L_1\cap
L_2^w}\sigma)^w\big).$$ 
Le lemme \ref{0.2.7} implique alors que $\r_{L_2}^G\pi$ est une somme
de repr\'esentations de $L_2$ qui sont des induites de \ecis\ de
sous-groupes de Levi standard de $L_2$ du type $^wL_1\cap L_2$.  
 Notons $W(L_1;L_2)_d$ l'ensemble des \'el\'ements $w\in W(L_1;L_2)$
tels que  
$^wL_1\cap L_2$ se transf\`ere (i.e. que $d$ divise le cardinal
de toutes ses sections). \'Ecrivons 
\begin{equation}\label{tttt}
\r_{L_2}^G\pi=\r_{L_2}^G\iii_{L_1}^G\sigma=
\end{equation}
$$
\sum_{w\in W(L_1,L_2)_d}\iii_{^wL_1\cap L_2}\big( ^w(\r_{L_1\cap
L_2^w}\sigma)\big)+\sum_{w\in W(L_1,L_2)\backslash
W(L_1,L_2)_d}\iii_{^wL_1\cap L_2}\big( ^w(\r_{L_1\cap
L_2^w}\sigma)\big)$$ 
D'abord, si $L_1$ ou $L_2$ ne se transf\`ere pas, alors
$W(L_1;L_2)_d=\O$ parce qu'un \sgls\ qui ne se transf\`ere pas ne peut
pas contenir un \sgls\ qui se transf\`ere.\\ 
 
C'est le cas en particulier si $L_1$ ne se transf\`ere pas
(i.e. $\pi\in {\mathcal B}_G\backslash \jl({\mathcal B}_{G'})$). On 
trouve donc dans ce cas que $\r_{L_2}^G\pi$ est la somme des induites
de \ecis\ de sous-groupes de Levi standard de $L_2$ qui ne se
transf\`erent pas. Cela veut dire que, si $L_2$ a $k$-blocs, l'image
de $\r_{L_2}^G : \rr(F)\to \otimes^k \rr(F)$ est incluse dans le
sous-espace  
$$\sum_{i=0}^{k-1} [\otimes^i\rr(F) \otimes \ii \otimes^{k-1-i}\rr(F)$$
de $\otimes^{k}\rr(F)]$. \\ 

En particulier, en appliquant la d\'efinition de $c$ on trouve que 
$$c(\pi)\in \rr(F)\otimes \ii+\ii\otimes\rr(F).$$
 Puisque $\pi$ est un \'el\'ement quelconque de ${\mathcal
B}_G\backslash \jl({\mathcal B}_{G'})$ cela implique que la
comultiplication ``passe au quotient'' dans $\rr(F)/\ii$. On note
$\bar{c}$ la comultiplication de $\rr(F)/\ii$ induite par $c$.\\ 

\'Etudions maintenant le cas o\`u $L_1$ se transf\`ere,
i.e. $\pi\in \jla({\mathcal B}_{G'})$. Par la formule \ref{tttt},
si $\pi$ est un  
\'el\'ement de $\jl({\mathcal B}_{G'})$, si $L_2$ a $k$ sections, dans
$\otimes^k[\rr(F)/\ii]$ on a~:  
$$(*)\ \ \ \ \ \ \ \ \ \ \ \r_{L_2}^G\pi=\sum_{w\in
W(L_1,L_2)_d}\iii_{^wL_1\cap L_2}\big( ^w(\r_{L_1\cap
L_2^w}\sigma)\big).$$ 

Si $L_2$ ne se transf\`ere pas, on a vu que $W(L_1,L_2)_d$ \'etait
vide et donc$$(***)\ \ \ \ \ \ \ \ \ \ \ \ \ \r_{L_2}^G\pi=0\
\text{dans}\  \otimes^{k}[\rr(F)/\ii].$$  

Si $L_2$ se transf\`ere, adoptons les notations suivantes~:

- soit $L'_1$ le \sgls\ de $G'$ correspondant \`a $L_1$, 

- soit $L'_2$ le \sgls\ de $G'$ correspondant \`a $L_2$, 

- soit $\pi'$ l'\'el\'ement de $Grot(G')$ correpondant \`a $\pi$~;
  autrement dit, $\pi'$ est l'induite de $L'_1$ \`a $G'$ de $\sigma'$
  o\`u $\sigma'$ est la \eci\ de $L'_1$ qui correspond \`a $\sigma$
  par la correspondance entre $L_1$ et $L'_1$. 

On a par le lemme \ref{geom}~:
$$(**)\ \ \ \ \ \ \ \ \ \ \ \ \r_{L'_2}^{G'}\pi'=\sum_{w\in
W(L'_1,L'_2)}\iii_{^wL'_1\cap L'_2}\big( ^w(\r_{L'_1\cap
L^{'w}_2}\sigma')\big).$$ 

Maintenant, il existe une inclusion standard $t$ de $W(L'_1;L'_2)$
dans $W(L_1;L_2)$. Elle est d\'efinie de la fa\c{c}on suivante. On
regarde les \'el\'ements de  
$W(L'_1;L'_2)$ comme des permutations de l'ensemble $\{1,2\ ...r\}$ et
les \'el\'ements de $W(L_1;L_2)$ comme des permutations de l'ensemble
$\{1,2\ ...n\}$ comme on l'a expliqu\'e au debut de cette
sous-section. Soit $\tau'$ une permutation dans $W(L'_1;L'_2)$. Pour
tout \'el\'ement $x$ de l'ensemble $\{1,2..n\}$ il existe un unique
couple $(a;b)\in \{1,2\ ...r\}\times \{0,1\ ...d-1\}$ tel qu'on ait   
$$x=ad-b.$$
On pose alors
$$\tau(x)=\tau'(a)d-b.$$
On obtient ainsi une permutation $\tau$ de l'ensemble $\{1,2\ ...n\}$
et on pose  
$$t(\tau')=\tau.$$
En fait, $t(\tau)$ est la permutation de l'ensemble $\{1,2\ ...n\}$ qui permute
les $r$ sous-ensembles $\{1,2\ ...d\},\ \{d+1,d+2\ ...2d\},\
...\{(r-1)d+1,(r-1)d+2\ ...n\}$ selon la permutation $\tau'$ tout en
laissant l'ordre inchang\'e \`a l'int\'erieur de chacun de ces
sous-ensembles. L'application clairement injective 
$$t : W(L'_1;L'_2)\to W(L_1;L_2)$$
ainsi d\'efinie a la propri\'et\'e que pour tout $w\in W(L'_1;L'_2)$,
les sous-groupes de Levi standard $^wL'_1\cap L'_2$ et $^{t(w)}L_1\cap
L_2$ se correspondent (et pareil pour $L'_1\cap L^{'w}_2$ et $L_1\cap
L^{t(w)}_2$). Ainsi, l'image par $t$ de $W(L'_1;L'_2)$ est incluse
dans $W(L_1;L_2)_d$. Le but du jeu est maintenant de montrer que
l'image par $t$ de $W(L'_1;L'_2)$ est en fait $W(L_1;L_2)_d$. Cela
donnera l'egalit\'e entre les membres droits des formules \ref{tttt}
et (**), gr\^ace au lemme \ref{0.2.7} b). 

On doit montrer que, si $w\in W(L_1;L_2)_d$ alors la permutation $w$
de l'ensemble $\{1,2\ ...n\}$ permute 
les sous-ensembles $\{1,2\ ...d\},\ \{d+1,d+2\ ...2d\},\
...\{(r-1)d+1,(r-1)d+2\ ...n\}$ en laissant l'ordre inchang\'e \`a
l'int\'erieur de chacun de ces sous-ensembles.  

Posons $w(1)=l$. Supposons que $l$ se trouve dans la $s$-i\`eme
section de $L_2$. Alors $l$ est le premier \'el\'ement de cette
section parce que, si $l'$ est un \'el\'ement de cette m\^eme section
de $L_2$ qui pr\'ec\`ede $l$, on doit avoir $w^{-1}(l')<w^{-1}(l)$
(voir d\'efinition de $W(L_1;L_2)$) ce qui est impossible car
$w^{-1}(l)=1$. Donc, $w(1)$ est un entier du type $\alpha d+1$,
$\alpha\in \{0,1\ ...r\}$ car il est le premier d'une section de $L_2$
et les cardinaux des sections de $L_2$ sont tous divisibles par $d$
puisque $L_2$ se transf\`ere.   

Montrons qu'on a $w(2)= \alpha d+2$, $w(3)=\alpha d+3$
...$w(d)=(\alpha+1)d$.  Montrons d'abord que les entiers $w(1)$,
$w(2)$ ...$w(d)$ se trouvent tous dans la m\^eme section de
$L_2$. Soit $\sigma$ la transposition $(1;i)$ pour un $i\in \{2,3\
...d\}$. On note toujours $\sigma$ la matrice de permutation qui lui
correspond. La matrice $\sigma$ se trouve dans le plus petit \sgls\ de
$G$ qui se transf\`ere (celui form\'e de matrices diagonales par blocs
de taille $d$ et qui est contenu par cons\'equence dans tout \sgls\ de
$G$ qui se transf\`ere). 
Donc $\sigma\in L_1\cap L_2^w$ puisque ce dernier est un \sgls\ qui se
transf\`ere. C'est \`a dire que $\sigma\in L_2^w$ soit $w\sigma
w^{-1}\in L_2$. Cela implique que la permutation  $w\sigma w^{-1}$
v\'erifie~: pour tout $x\in \{1,2\ ...n\}$, $w\sigma w^{-1}(x)$ et $x$
se trouvent dans la m\^eme section de $L_2$. 
En appliquant cela \`a $x=w(1)$ on trouve que $w(i)$ et $w(1)$ se
trouvent dans la m\^eme section de $L_2$. C'est ce qu'on voulait
d\'emontrer.

Maintenant prouvons que $w(2)= \alpha d+2$, $w(3)=\alpha d+3$ ... et
$w(d)=(\alpha+1)d$. Si, par exemple, $w(2)\neq  \alpha d+2$, alors on
a $w(1)=\alpha d+1<\alpha d+2<w(2)$. Ces trois entiers \'etant dans la
m\^eme section de $L_2$ (puisque $w(1)$ et $w(2)$ le sont), on a
$w^{-1}(w(1))<w^{-1}(\alpha d+2)<w^{-1}(w(2))$ (voir encore une fois
la d\'efinition de $W(L_1;L_2)$) ce qui est impossible car il n'y a
aucun entier entre 1 et 2. 
Ainsi de suite pour $w(3)$ et les autres. On prouve ainsi que $w$
envoie $\{1,2\ ...d\}$ sur $\{\alpha d+1,\alpha d+2\ ...(\alpha+1)d\}$
sans changer l'ordre.  

On continue de la fa\c{c}on suivante~: on prouve que $w(d+1)$ est un
entier du type $\beta d+1$. En effet, s'il n'est pas le premier dans
la section de $L_2$ qui le contient, les seuls entiers qui puissent le
pr\'ec\'eder dans cette section sont $w(1)$, $w(2)$ ... et
$w(d)$. Mais ces entiers sont cons\'ecutifs donc si un le pr\'ec\`ede,
tous le pr\'ec\`edent et on a alors 
$w(d+1)=w(d)+1=(\alpha +1)d+1$. Apr\`es cette remarque la
d\'emonstration se d\'eroule comme pour $\{1,2\ ...d\}$, jusqu'au
r\'esultat final.\\ 

On a trouv\'e donc que $t$ r\'ealise une bijection~:
$$t : W(L'_1;L'_2)\simeq W(L_1;L_2)_d.$$
Si $L_2$ a $k$ sections, notons $\jlh^{\otimes^k}$ l'isomorphisme d'anneaux 
$$\otimes^k\rr(D) \simeq \otimes^k [\rr(F)/\ii]$$
 induit par l'isomorphisme d'anneaux 
$$\jlh : \rr(D)\simeq \rr(F)/\ii.$$
 Alors l'image de $\r_{L_2}^G(\pi)$ dans $\otimes^k[\rr(F)/\ii]$ est
\'egale \`a $\jlh^{\otimes^k}(r_{L'_2}^{G'}(\pi'))$. Cela se voit en comparant
les formules (*) et (**)~: pour tout $w\in W(L'_1;L'_2)$, 
$$\iii_{^wL'_1\cap L'_2}\big( ^w(\r_{L'_1\cap L^{'w}_2}\sigma)\big)\
\text{et}\  \iii_{^{t(w)}L_1\cap L_2}\big( ^{t(w)}(\r_{L_1\cap
L_2^{t(w)}}\sigma')\big)$$ 
 se correspondent par $\jlh^{\otimes^k}$ et $t$ est une bijection de
$W(L'_1;L'_2)$ sur $W(L_1;L_2)_d$. Maintenant revenons aux
d\'efinitions de la comultiplication et de l'involution. Notons
${\mathcal L}$ l'ensemble de sous-groupes de Levi standard de $G$ et
${\mathcal L}_{max}$ son sous-ensemble form\'e de sous-groupes de Levi
standard maximaux. De m\^eme, notons ${\mathcal L}'$ l'ensemble de
sous-groupes de Levi standard de $G'$ et ${\mathcal L}'_{max}$ son
sous-ensemble form\'e de sous-groupes de Levi standard
maximaux. Supposons qu'on est dans la situation $L_1$ se transf\`ere
en $L'_1$, $\pi=\iii_{L_1}^G\sigma$ o\`u $\sigma$ est une \eci\ de
$L_1$, $\sigma'$ est la \eci\ de $L'_1$ qui lui correspond et
$\pi'=\iii_{L'_1}^{G'}\sigma'$.  

On a
$$c(\pi)=\sum_{L\in {\mathcal L}_{max}}\r_L^G(\pi)\in \rr(F)\otimes\rr(F).$$
Cette somme se d\'ecompose en deux sommes~: une sur les sous-groupes
de Levi standard maximaux de $G$ qui se transf\`erent et l'autre sur
les sous-groupes de Levi standard maximaux de $G$ qui ne se
transf\`erent pas. L'image (bien d\'efinie, on l'a vu) de $c(\pi)$
dans $\rr(F)/\ii\otimes \rr(F)/\ii$ est la m\^eme que l'image de la
somme sur les sous-groupes de Levi standard maximaux de $G$ qui se
transf\`erent par la remarque (***) plus haut. 
D'autre part on a  
$$c'(\pi')=\sum_{L'\in {\mathcal L}'_{max}}\r_{L'}^{G'}(\pi')\in \rr(D)\otimes\rr(D).$$
Les sous-groupes de Levi standard maximaux de $G$ qui se transf\`erent
se correspondent biunivoquement avec les sous-groupes de Levi standard
maximaux de $G'$. Si $L$ est un tel sous-groupe de $G$ et $L'$ est le
sous-groupe de $G'$ qui lui correspond, alors on a vu que l'image de
$\r_L^G\pi$ dans $\rr(F)/\ii\otimes \rr(F)/\ii$ est \'egale \`a
$\jlh^{\otimes^2}(\r_{L'}^{G'}\pi')$. Donc, finalement, l'image $\bar{c}(\pi)$
de $c(\pi)$ dans $\rr(F)/\ii\otimes \rr(F)/\ii$ est \'egale \`a
$\jlh^{\otimes^2}(c(\pi'))$.\qed

\subsection{Une pr\'ecision}\label{precision}

Pour un r\'esultat plus pr\'ecis (prop. \ref{corresp}), 
il nous faut rappeler les travaux de
[Ta]. Soit $G'=GL_r(D)$ une forme int\'erieure de $G$. Si $\rho'$
est une \cusp\ de $G'$, alors il existe une repr\'esentation $\rho$
essentiellement de carr\'e int\'egrable de $G$ qui correspond \`a
$\rho'$ par Jacquet-Langlands. Cette repr\'esentation $\rho$ se trouve
au-dessus d'un 
segment de Zelevinski de longueur $k$. On pose alors
$s(\rho')=k$. 
Notons $\nu$ la norme r\'eduite sur $G'$. Si $\rho'$ est une
repr\'esentation cuspidale de $G'$ on pose
$\nu_{\rho'}=\nu^{s(\rho')}$. Si $\sigma'$
est une repr\'esentation essentiellement de carr\'e int\'egrable de
$G'$, alors il existe une repr\'esentation cuspidale $\rho'$ de
$GL_t(D)$, o\`u $t$ divise $r$, telle que $\sigma'$ soit l'unique
quotient de l'induite \`a $G'$ de la repr\'esentation $\rho'\otimes
\nu_{\rho'}\rho'\otimes
\nu^{2}_{\rho'}\rho'...\otimes
\nu^{(\frac{r}{t}-1)}_{\rho'}\rho'$ vue comme repr\'esentation du groupe
form\'e des matrices par blocs diagonaux de taille $t$. On note alors
$\sigma'$ \`a la Zelevinski, $<[\rho',
\nu^{(\frac{r}{t}-1)}_{\rho'}\rho']>^t$.
On appelle la
 repr\'esentation $\rho'\otimes
\nu_{\rho'}\rho'\otimes
\nu^{2}_{\rho'}\rho'...\otimes
\nu^{(\frac{r}{t}-1)}_{\rho'}\rho'$
segment de Tadi$\check{\rm c}$ (de Zelevinski
quand c'est $GL_n$ 
lui-m\^eme) et nous l'assimilons \`a l'ensemble
$\{\rho',\nu_{\rho'}\rho',\nu^{2}_{\rho'}\rho'...
\nu^{(\frac{r}{t}-1)}_{\rho'}\rho'\}$ quand il s'agit de faire des
r\'eunions et intersections de segments. 
 Posons par convention que
l'ensemble vide est aussi un segment. On dit que
deux segments $\Delta_1$ et $\Delta_2$ de
Tadi$\check{\rm c}$ sont {\it li\'es} si  leur r\'eunion en tant qu'ensembles est
 un segment diff\'erent \`a la fois de $\Delta_1$ et $\Delta_2$.
 On note ce segment
$\Delta_1\cup\Delta_2$. L'intersection est toujours un segment qu'on note
$\Delta_1\cap\Delta_2$. 

Soient $L'$ un \sgls\ de $G'$ qu'on identifie \`a un
produit $GL_{n_1}(D)\times GL_{n_2}(D)\times ...\times GL_{n_k}(D)$,
et $\sigma'= \sigma'_1\otimes \sigma'_2\otimes...\otimes \sigma'_k$
une \eci\ de $L'$.

Pour chaque $i$, soit $\Delta_i$ le segment correspondant \`a
$\sigma'_i$. Si parmi les $\Delta_i$ il y en a deux qui sont li\'es,
$\Delta_{i_1}$ et $\Delta_{i_2}$, 
disons que nous faisons une {\it op\'eration \'el\'ementaire} sur l'ensemble
des $\Delta_i$ si on remplace les deux segments $\Delta_{i_1}$ et
$\Delta_{i_2}$ par  les deux autres segments
$\Delta_{i_1}\cap\Delta_{i_2}$ $\Delta_{i_1} \cup
\Delta_{i_2}$. Au nouvel ensemble de segments ainsi obtenu correspond
une \eci\ $\tau'$ d'un \sgls\ de $G'$ (modulo permutation des blocs plus
pr\'ecis\'ement). Nous disons alors que $\tau'$ a \'et\'e obtenue
de $\sigma'$ par
une op\'eration \'el\'ementaire. Si $\tau'$ et $\sigma'$ sont des
\ecis\ de sous-groupes de Levi standard de $G'$, on dit que $\tau'$
est inf\'erieure \`a $\sigma'$ si elle peut s'obtenir de
$\sigma'$ par un nombre fini d'op\'erations \'el\'ementaires. 
Ceci induit un ordre partiel sur l'ensemble des \ecis\ des \sglss\
modulo conjugaison. Du coup, nous avons une relation d'ordre \'evidente
associ\'ee sur ${\mathcal B}_{G'}$. 

Aussi, pour
toute \eci\ $\sigma'$ de $G'$ il existe un unique nombre r\'eel
$e(\sigma')$ tel que la repr\'esentation $\nu^{e(\sigma')}\sigma'$
soit unitaire.
Soit  $\sigma'= \sigma'_1\otimes \sigma'_2\otimes...\otimes \sigma'_k$
tel que 
$e(\sigma'_i)\geq e(\sigma'_{i+1})$ pour tout $1\leq i\leq k-1$. 
Alors l'induite de $\sigma'$ \`a $G'$ a un unique
quotient irr\'eductible. On le note $Lg(\sigma')$. La notation se
justifie par le fait que, regroupant les $\sigma'_i$ par des paquets
maximaux qui ont le m\^eme coefficient $e(\ )$, on peut induire $\sigma'$
\`a un sous-groupe parabolique tel que l'induite soit irr\'eductible
et essentiellement temp\'er\'ee, et $Lg(\sigma')$ soit le quotient de
Langlands de cette repr\'esentation. Pour toute repr\'esentation lisse
irr\'eductible $\pi'$ de $G'$ il existe un $\sigma'$ comme plus haut
tel que $\pi'$ soit $Lg(\sigma')$. La repr\'esentation $\sigma'$ est
unique modulo l'op\'eration qui consiste \`a permuter entre eux des
$\sigma'_i$ qui ont le m\^eme coefficient $e$. On voit que
l'application 
$$Q : {\mathcal B}_{G'}\to Irr(G')$$
donn\'ee par ${\bf
i}_{L'}^{G'}(\sigma')\mapsto Lg(\sigma')$ est bien d\'efinie et
bijective. La relation d'ordre sur  ${\mathcal B}_{G'}$ induit ainsi
une relation d'ordre sur $Irr(G')$. On \'etudie plus en d\'etail la
relation entre ces deux bases remarquables de $Grot(G')$ que sont
${\mathcal B}_{G'}$ et $Irr(G')$.

On utilise le r\'esultat fondamental suivant
 ([Ta], g\'en\'eralisant [Ze])~:
\begin{prop}\label{tadic}
Soit $\sigma'$ une \eci\ d'un \sgls\ de $G'$. Dans l'\'ecriture de
l'image dans $Grot(G')$ de l'induite de $\sigma'$ 
\`a $G'$ comme somme de classes de repr\'esentations irr\'eductibles, 

- $Lg(\sigma')$ apparait avec multiplicit\'e exactement 1

- une repr\'esentation irr\'eductible $\pi'$ appara\^\i t avec une
  multiplicit\'e non nulle si et seulement si $\pi'\leq \sigma'$.
\end{prop}

Si $L'$ est un \sgls\ de $G'$ et $\rho'$ est une
repr\'esentation cuspidale de $L'$, notons $Grot(G')_{(L',\rho')}$ le
sous-groupe de $Grot(G')$ engendr\'e par les images des repr\'esentations
irr\'eductibles qui sont des sous-quotients de ${\bf
i}_{L'}^{G'}(\rho')$. Soit $X$ l'ensemble de classes modulo
permutation par blocs de couples $(L',\rho')$. Il est clair que,
si $x,y\in X$ et $x\neq y$ alors $Grot(G')_x$ et $Grot(G')_y$ sont 
disjoints
et que $Grot(G')$ est
la somme directe des $Grot(G')_{x}$ sur tous les $x\in X$. Pour tout
$x\in X$, $Grot(G')_{x}$ est un $\z$-module de dimension finie. Il
admet deux bases naturelles~: ${\mathcal B}_{G'}\cap Grot(G')_{x}$ et
$Irr(G')\cap Grot(G')_{x}$. Ordonnons la base ${\mathcal
B}_{G'}\cap Grot(G')_{x}$ dans un ordre croissant pour la relation
d\'efinie plus haut (i.e. chaque fois qu'un \'el\'ement est
inf\'erieur \`a un autre, il le pr\'ec\`ede). Ordonnons la base
$Irr(G')\cap Grot(G')_{x}$ dans l'ordre (croissant) correspondant via
l'aplication $Q$. La matrice de passage entre les deux bases est de
taille finie, 
et la proposition \ref{tadic}   
implique qu'elle est triangulaire sup\'erieure unipotente. C'est donc
pareil pour la matrice de passage dans l'autre sens, ce qui prouve
le r\'esultat suivant, dual en quelque sorte du pr\'ec\'edent~:

\begin{prop}\label{eu}
Soit $\pi'=Lg(\sigma')$ ($\sigma'$ \eci\ d'un \sgls\ $L'$) une
repr\'esentation lisse irr\'eductible de 
$G'$. Dans $Grot(G')$ on a 
$$\pi'={\bf i}_{L'}^{G'}(\sigma')+\sum_{\sigma'_i<\sigma'}
a_i{\bf i}_{L'_i}^{G'}(\sigma'_i)$$
o\`u les $a_i$ sont des entiers et, pour tout $i$, $L'_i$ est un
\sgls\ et $\sigma'_i$ est une \eci\ de $L'_i$.
\end{prop}

Remarquons au passage que, si $x$ et $y$ appartenant \`a $X$ sont
distincts, alors aucun \'el\'ement de la base ${\mathcal
B}_{G'}\cap Grot(G')_{x}$ n'est en relation d'ordre avec un
\'el\'ement de la base ${\mathcal
B}_{G'}\cap Grot(G')_{y}$. C'est une cons\'equence du fait qu'en
faisant une op\'eration \'el\'ementaire on reste ``au-dessus'' d'une
m\^eme classe cuspidale $(L',\rho')$.

\begin{lemme}\label{ordre}
L'application $\jl$ de ${\mathcal
B}_{G'}$ dans ${\mathcal
B}_{G}$ est croissante.
\end{lemme}

{\bf D\'emonstration.} L'ensemble des segments de Tadi$\check{\rm c}$
de $G'$ est en 
bijection avec l'ensemble des \ecis\ de $G'$, et l'ensemble des
segments de Zelevinski de $G$ est en
bijection avec l'ensemble des \ecis\ de $G$. \`A travers la
correspondance de Jacquet-Langlands il existe donc une
bijection naturelle de l'ensemble 
des segments de Tadi$\check{\rm c}$ de $G'$ dans l'ensemble de segments de
Zelevinski de $G$ (attention, elle ne co\"\i ncide pas avec la
restriction de 
l'application $\jlh:\rr(D)\to\rr(F)$).

Il est
facile \`a voir
qu'\`a des segments de Tadi$\check{\rm c}$ $\Delta_1$ et $\Delta_2$
li\'es sur $G'$ 
correspondent par cette bijection des segments de 
Zelevinski
$\Gamma_1$ et $\Gamma_2$
li\'es sur $G$ (la r\'eciproque \'etant 
fausse). Aussi, \`a  $\Delta_1\cap\Delta_2$ correspond
$\Gamma_1\cap\Gamma_2$ et \`a  $\Delta_1\cup\Delta_2$ correspond
$\Gamma_1\cup\Gamma_2$. Ceci montre que notre injection de  ${\mathcal
B}_{G'}$ dans ${\mathcal B}_G$ pr\'eserve
l'ordre d\'efini plus haut.\qed

On pourrait penser que, si $a,b\in {\mathcal
B}_{G'}$, si $\jl(b)<\jl(a)$, alors $b<a$.
Ceci n'est pas vrai, comme le montre le contre-exemple suivant qui m'a
\'et\'e fourni par Andrei
Moroianu:\\
\ \\
{\bf Fait (troublant).} Consid\'erons les repr\'esentations
$$\sigma_1=ind(<[\rho;\nu^5\rho]>^t\otimes <[\nu\rho;\nu^8\rho]>^t\otimes
<[\nu^3\rho;\nu^6\rho]>^t\otimes <[\nu^4\rho;\nu^9\rho]>^t)$$
et
$$\sigma_2=ind(<[\rho;\nu^9\rho]>^t\otimes <[\nu\rho;\nu^6\rho]>^t\otimes
<[\nu^3\rho;\nu^8\rho]>^t\otimes <[\nu^4\rho;\nu^5\rho]>^t)$$
de $G=GL_{24}(F)$, o\`u $\rho$ est un caract\`ere de $F^*$. 
$\sigma_1$ et $\sigma_2$ sont des \'el\'ements de ${\mathcal B}_G$ et 
on a $\sigma_2<\sigma_1$. Soit $D$ une alg\`ebre de quaternions sur
$F$. Soit $\rho'$ la repr\'esentation de $D^*$ qui correspond \`a
$<[\rho;\nu\rho]>^t$ par Jacquet-Langlands. On a $\nu_{\rho'}=\nu^2$.
On consid\`ere maintenant les repr\'esentations
$$a=ind(<[\rho';\nu_{\rho'}^2\rho']>^t\otimes
<[\nu\rho';\nu_{\rho'}^3\nu\rho']>^t\otimes 
<[\nu^3\rho';\nu_{\rho'}\nu^3\rho']>^t\otimes 
<[\nu^4\rho';\nu_{\rho'}^2\nu^4\rho']>^t)$$
et 
$$b=ind(<[\rho';\nu_{\rho'}^4\rho']>^t\otimes
<[\nu\rho';\nu_{\rho'}^2\nu\rho']>^t\otimes 
<[\nu^3\rho';\nu_{\rho'}^2\nu^3\rho']>^t\otimes 
<[\nu^4\rho';\nu^4\rho']>^t)$$

de $G'=GL_{12}(D)$. Alors $a,b\in  {\mathcal
B}_{G'}$, $\jl(a)=\sigma_1$, $\jl(b)=\sigma_2$, mais il est clair
qu'on n'a pas $b<a$, (seuls les segments
$[\rho';\nu_{\rho'}^2\rho']$ et $[\nu^4\rho';\nu_{\rho'}^2\nu^4\rho']$
dans la formule de $a$ sont li\'es).\\

Donnons maintenant le r\'esultat qui pr\'ecise un peu plus la
correspondance $\jl$~:

\begin{prop}\label{corresp}
Soit $\pi'$ une repr\'esentation lisse irr\'eductible de
$G'$. Posons $\pi'=Lg(\sigma')$, o\`u $\sigma'$ est une \eci\ d'un
\sgls\ $L'$ de $G'$. Soient $L$ le \sgls\ de $G$ qui correspond \`a
$L'$, et $\sigma={\bf C}^{-1}(\sigma')$ la \eci\ de $L$ qui correspond \`a
$\sigma'$. Posons $\pi=Lg(\sigma)$.
Alors on a
$$\jl(\pi')=\pi+\sum_{\sigma_i<\sigma}a_i{\bf i}_{L_i}^{G}(\sigma_i)
=\pi+\sum_{\pi_j<\pi}b_j\pi_j$$
o\`u les $a_i$ et les $b_j$ sont des entiers et, pour tout $i$, $L_i$ est un
\sgls\ de $G$, $\sigma_i$ sont des \eci\ de $L_i$. et $\pi_j$ sont des
repr\'esentations lisses irr\'eductibles de $G$. 
\end{prop}

{\bf D\'emonstration.} \'Ecrivons 
$$\pi'=
{\bf i}_{L'}^{G'}(\sigma')+\sum_{\sigma'_i<\sigma'}a_i{\bf i}_{L'_i}^{G'}
(\sigma'_i)$$
comme dans la prop. \ref{eu}. La d\'efinition de $\jl$ donne~: 
\begin{equation}\label{scris}
\jl(\pi')={\bf i}_{L}^{G}(\sigma)+\sum_{\sigma_i}a_i{\bf
i}_{L_i}^{G}(\sigma_i)
\end{equation}
o\`u $L_i$ correspond \`a $L'_i$ et $\sigma_i$ correspond \`a
$\sigma'_i$. On a $\sigma_i<\sigma$ par le lemme \ref{ordre}.

Appliquons la prop. \ref{eu} \`a $\pi$~:
$$\pi={\bf i}_{L}^{G}(\sigma)+\sum_{\sigma_j<\sigma}a'_j{\bf
i}_{M_j}^{G}(\sigma_j)$$
et substituons ${\bf i}_{L}^{G}(\sigma)$ dans l'\'equation
\ref{scris}. On obtient
$$\jl(\pi')=\pi-[\sum_{\sigma_j<\sigma}a'_j{\bf
i}_{M_j}^{G}(\sigma_j)]+
\sum_{\sigma_i<\sigma}a_i{\bf i}_{L_i}^{G}(\sigma_i).$$
Cela prouve la premi\`ere \'egalit\'e. La deuxi\`eme est \'evidente.
\qed 

\subsection{Une conjecture}

\`A ce stade, on peut se poser la question, si $\pi'$ et $\pi$ sont
comme plus haut, si on a la relation entre caract\`eres 
$$\chi_{\pi}(g)=(-1)^{n-r}\chi_{\pi'}(g')$$
pour tout $g$ et $g'$ qui se correspondent. Nous ne pourrons pas d\'emontrer
que $\jl(\pi')=\pi$, puisqu'on sait que $\jl$ n'envoie pas toute 
repr\'esentation irr\'eductible sur une repr\'esentation irr\'eductible,
mais il se peut que $\pi'$ et $\pi$ se correspondent toutefois. Dans le 
contre-exemple du deuxi\`eme commentaire \`a la section \ref{grot} nous 
d\'efinissons $\pi'$ irr\'eductible
telle que $\jl(\pi')$ ne soit pas une repr\'esentaion
irr\'eductible, mais, dans ce cas par exemple,
$\pi'$ correspond \`a $\pi$ quand-m\^eme ! (voir prop. \ref{2} plus bas).

\subsubsection{La conjecture}
Ceci nous motive \`a poser une conjecture. On va dire que $\pi\in
Grot(G)$ et $\pi'\in Grot(G')$ se correspondent si on a la relation
entre les caract\`eres 
\begin{equation}\label{caractere}
\chi_{\pi}(g)=(-1)^{n-r}\chi_{\pi'}(g')
\end{equation}
pour tous $g\in G$ et $g'\in G'$ qui se correspondent.\\
\ \\
{\bf Conjecture.} {\it Soit $\pi'=Lg(\sigma')$ une repr\'esentation
lisse irr\'eductible 
de $G'$. Alors $\pi'$ correspond \`a $\pi=Lg(\sigma)$, o\`u
$\sigma={\bf C}^{-1}(\sigma')$.}\\ 

Nous ne savons pas prouver cette conjecture en g\'en\'eral.\\

\subsubsection{Possible application} 

Tout d'abord, que la conjecture soit vraie ou pas, on pose
$M(\pi')=\pi$. Il est \'evident que
$\pi'\mapsto M(\pi')$ induit par lin\'earit\'e un
morphisme injectif de groupes $M:Grot(G')\to Grot(G)$, qui respecte
l'irr\'eductibilit\'e ($M$ n'est autre que $Q^{-1}\jl Q$). 
Nous noterons abusivement par la suite
avec la m\^eme lettre $M$ ce morphisme pour tous les entiers $r$
confondus.

\begin{prop}\label{applicconj}
Supposons que la conjecture soit 
vraie.
Soient $\pi'_1\in Irr(GL_{r_1}(D))$ et $\pi'_2\in
Irr(GL_{r_2}(D))$. Si
l'induite \`a $GL_{d(r_1+r_2)}(F)$ de $M(\pi'_1)\otimes M(\pi'_2)$ 
est irr\'eductible, alors l'induite \`a $GL_{(r_1+r_2)}(D)$ de
$\pi'_1\otimes \pi'_2$ est 
irr\'eductible (se g\'en\'eralise imm\'ediatement \`a plusieurs
repr\'esentations).
\end{prop}

{\bf D\'emonstration.} Posons $\pi'_1=Lg(\sigma'_1)$,
$\pi'_2=Lg(\sigma'_2)$, $\pi_1=Lg(\sigma_1)$ et $\pi_2=Lg(\sigma_2)$
($\sigma_i$ et $\sigma'_i$ se correspondent). Si
l'induite \`a $GL_{d(r_1+r_2)}(F)$ de $M(\pi'_1)\otimes M(\pi'_2)$ 
est irr\'eductible, elle est \'egale \`a
$Lg(\sigma_1\otimes\sigma_2)$ ([Ta]). Donc $Lg(\sigma_1\otimes\sigma_2)$ 
correspond \`a la fois \`a l'induite de $\pi'_1\otimes\pi'_2$ et et
\`a $Lg(\sigma'_1\otimes\sigma'_2)$ (par la conjecture suppos\'ee
vraie). Ces deux repr\'esentations de $GL_{(r_1+r_2)}(D)$ ont donc des
caract\`eres \'egaux, elles sont donc dans la m\^eme classe.\qed
\ \\
\subsubsection{Cas (tr\`es) particuliers}

\begin{prop}\label{2}
Soient $\sigma'_1$ et $\sigma'_2$ deux
repr\'esentations \eci\ de $GL_{r_1}(D)$ et $GL_{r_2}(D)$
respectivement, $r_1+r_2=r$. Alors tout sous-quotient 
irr\'eductible de
l'induite \`a $G'$ de $\sigma'_1\otimes \sigma'_2$ v\'erifie 
la conjecture.
\end{prop}

{\bf D\'emonstration.} Soient 
$\sigma_1=\jl (\sigma'_1)$ et $\sigma_2=\jl (\sigma'_2)$.
On a 
$$\jl(ind\, \sigma'_1\otimes \sigma'_2)=ind\, \sigma_1\otimes \sigma_2$$
donc ces deux repr\'esentations se correspondent. 

Si les segments associ\'es \`a $\sigma_1$ et $\sigma_2$ 
ne sont pas li\'es alors les segments associ\'es \`a 
$\sigma'_1$ et $\sigma'_2$ ne sont pas li\'es non plus, les deux induites
sont irr\'eductibles et la proposition est d\'emontr\'ee. 

Si 
les segments associ\'es \`a $\sigma_1$ et $\sigma_2$ 
 sont li\'es, alors 
$$ind\, \sigma_1\otimes \sigma_2=Lg(\sigma_1\otimes \sigma_2)+ind\, \tau$$
o\`u $\tau$ est la seule r.e.c.i. inf\'erieure \`a $\sigma_1\otimes \sigma_2$.
Par ailleurs $\tau$ est minimale et donc son induite est irr\'eductible. On a
donc
$$ind\, \sigma_1\otimes \sigma_2=Lg(\sigma_1\otimes \sigma_2)+Lg(\tau)$$
Maintenant on a deux possibilit\'es~: 

- ou bien les segments associ\'es 
\`a $\sigma'_1$ et $\sigma'_2$ sont li\'es, ce qui arrive si et seulement
si le sous-groupe de d\'efinition de $\tau$ se transf\`ere, et alors 
$$ind\, \sigma'_1\otimes \sigma'_2=Lg(\sigma'_1\otimes \sigma'_2)+ind\, \tau'
.$$
La r.e.c.i. $\tau'$ est minimale et on \'ecrit 
$$ind\, \sigma'_1\otimes \sigma'_2=Lg(\sigma'_1\otimes \sigma'_2)+Lg(\tau').$$

On a $\tau=\jl(\tau')$ donc leurs induites se correspondent, ce qui montre
que $Lg(\tau')$ correspond \`a $Lg(\tau)$. Donc aussi, 
n\'ecessairement, $Lg(\sigma'_1\otimes \sigma'_2)$ et 
$Lg(\sigma'_1\otimes \sigma'_2)$ se correspondent.

- ou bien les segments asoci\'es 
\`a $\sigma'_1$ et $\sigma'_2$ ne sont pas li\'es, i.e. 
$$ind\, \sigma'_1\otimes \sigma'_2=Lg(\sigma'_1\otimes \sigma'_2)$$
est irr\'eductible. Mais, si les segments asoci\'es 
\`a $\sigma'_1$ et $\sigma'_2$ ne sont pas li\'es cela veut dire que le
groupe de d\'efinition de $\tau$ ne se transf\`ere pas, donc l'induite
\`a $G$ de $\tau$ est $G'$-nulle. Donc $Lg(\sigma'_1\otimes \sigma'_2)$
correspond \`a $Lg(\sigma_1\otimes \sigma_2)$.\qed

\begin{cor} 
Si $r=2$ la conjecture est v\'erifi\'ee.
\end{cor}

{\bf D\'emonstration.} Une  repr\'esentation de $GL_2(D)$ est ou bien 
cuspidale, ou bien un sous-quotient  comme dans la proposition.\qed

\begin{prop}\label{3}
Soit $S'=\{\rho',\nu_{\rho'}\rho'...\nu^k_{\rho'}\rho'\}$ un segment de $G'$.
Alors tout sous-quotient irr\'eductible 
de l'induite de $S$ \`a $G'$ v\'erifie 
la conjecture.
\end{prop}

{\bf D\'emonstration.} Soit $\pi'$ un tel sous-quotient. 
Notons $L'$ le sous-groupe de Levi du segment $S'$,
$L$ le sous-groupe de Levi de de $G'$ qui lui correspond, et
$S$ la \eci\ de $L$ qui correspond \`a $\rho'$. Il arrive que, dans
ce cas particulier, toutes 
les \ecis\ inf\'erieures \`a $S$ (du c\^ot\'e $G$) sont d\'efinies
sur des sous-groupes de Levi qui se transf\`erent, et correspondent
donc biunivoquement aux \ecis\ inf\'erieures \`a $S'$. 
De plus, on sait que, dans l'induite d'une telle \eci\, 
tous les sous-quotients apparaissent avec multiplicit\'e 
exactement $1$. On peut donc faire une d\'emonstration par r\'ecurrence
sur le nombre de repr\'esentations de $G'$ inf\'erieures \`a $\pi'$.

Si $\pi'$ est minimale, elle est une \eci\
et le r\'esultat d\'ecoule de la prop. \ref{2}. 

Supposons que,
pour tout $\tau'<\pi'$, $\tau'$ correspond \`a $M(\tau')$.
Si
$\pi'=Lg(\sigma')$ et $\sigma'$ correspond \`a $\sigma$,
alors les induites de $\sigma$ et $\sigma'$ se correspondent.
On en d\'eduit que
$\pi'+\sum_{\tau'<\pi'}\tau'$ correspond \`a $M(\pi')
+\sum_{\tau'<\pi'}M(\tau')$. D'o\`u le r\'esultat en appliquant l'hypoth\`ese
de r\'ecurrence.\qed

Un cas int\'eressant sont les repr\'esentations g\'en\'eriques de
$G$. On attend d'une repr\'esentation automorphe cuspidale de $GL_n$
que toutes ses composantes locales soient temp\'er\'ees aux
places non archim\'ediennes, mais pour l'instant nous savons du moins
qu'elles sont toujours g\'en\'eriques ([Sh], pp. 190).
Par [Ze], nous savons aussi qu'une repr\'esentation g\'en\'erique $\pi$
de $G$ est
toujours de la forme $ind_P^G(\sigma)$, o\`u $P$ est un sous-groupe
parabolique standard de $G$ et $\sigma$ est une r.e.c.i. du \sgls\ de
$P$ 
telle
que l'induite soit irr\'eductible. Nous avons la 
\begin{prop}\label{generiques}
Si $P$ ne se transf\`ere pas, alors le caract\`ere de $\pi$ est
$G'$-nul.
Si $P$ se transf\`ere, alors $\pi$ v\'erifie la conjecture.
\end{prop}

{\bf D\'emonstration.} La premi\`ere assertion est \'evidente, car
$\pi$ se trouve dans $\ss$ par d\'efinition m\^eme. La deuxi\`eme
affirmation n'est pas difficile~: si $P'$ est le sous-groupe
parabolique standard de $G'$ qui correspond \`a $P$, et $\sigma'$ est
la r.e.c.i. du \sgls\ de $P'$ qui correspond \`a $\sigma$ par
Jacquet-Langlands, alors les crit\`eres des r\'eductibilit\'e de [Ze]
et [Ta] pour les induites des r.e.c.i. impliquent que,
$ind_P^G(\sigma)$ \'etant irr\'eductible, 
$ind_{P'}^{G'}(\sigma')$ l'est aussi. D'o\`u le r\'esultat.\qed

Le transfert des composantes locales des
repr\'esentations automorphes cuspidales de $GL_n$ est donc simple. 
Le probl\`eme des {\it certaines} composantes locales des
repr\'esentations du spectre r\'esiduel 
est \'etudi\'e dans la sous-section suivante.

\subsection{L'involution}\label{inv} 

A. Zelevinski a d\'efini une involution de l'alg\`ebre de Hopf
$\rr(F)$ (prop. 9.12, [Ze]). Dans [Au], A.-M. Aubert d\'efinit une
involution du groupe de Grothendieck d'un groupe r\'eductif
quelconque. Elle montre \'egalement que cette involution envoie une
repr\'esentation irr\'eductible sur une repr\'esentation
irr\'eductible au signe pr\`es et que cette involution co\"\i ncide
au signe pr\`es avec celle de Zelevinski dans le cas du groupe
lin\'eaire. Posons $G=GL_n(F)$. Sur
$Grot(G)$, l'involution de [Au] est d\'efinie de la fa\c{c}on
suivante~: si $\mathcal L$ est l'ensemble de sous-groupes de Levi
standard de 
$G$, si $\pi$ est un \'el\'ement de $Grot(G)$, alors  
$$i(\pi)=\sum_{L\in {\mathcal
L}}(-1)^{|L|}\iii_L^{G}\r_{L}^{G}(\pi).$$ 
Cette involution ne s'\'etend pas bien \`a $\rr(F)$ comme celle de
Zelevinski, \`a cause d'un signe qui entrave la lin\'earit\'e ([Au],
th. 2.3). Nous alons \'etudier le comportement de $i$ par rapport au
morphisme $\overline{\jl}$ d\'efini au th. \ref{jlr}, c). 

Notons $i'$ involution de A.-M. Aubert pour $G'=GL_r(D)$ (formule
analogue \`a celle pour $G$).
\begin{theo}\label{involution}
L'involution $i$ de $Grot(G)$ induit une involution bien d\'efinie
$\bar i$ de
$Grot(G)/\ss$. L'isomorphisme $\overline{\jl}$ transforme l'involution
$i'$ en l'involution $(-1)^{n-r}\bar i$.
\end{theo}

{\bf D\'emonstration.} Soit $\pi\in \ss$. On veut montrer que $i(\pi)\in
\ss$. Or, on a vu que $\ii$ est stable par induction et par
restriction. Il suffit donc d'utiliser la formule
$$i(\pi)=\sum_{L\in {\mathcal L}}(-1)^{|L|}\iii_L^G\r_L^G(\pi),$$
car $\ss=Grot(G)\cap\ii$.
 
Ainsi, l'involution $i$ ``passe au quotient'', et $\bar i$ est bien
d\'efinie (involution, \'evidemment).

On montre que $\overline{\jl}(i')=(-1)^{n-r}\bar i$. Si $\pi\in
Grot(G')$, on veut prouver que
$$(-1)^{n-r} i (\jl(\pi))=\jl(i'(\pi)).$$

On utilise encore $\jlh$. En effet, il suffit de montrer que 
$$(-1)^{n-r} i (\jlh(\pi))=\jlh(i'(\pi)).$$

On a 
$$i'(\pi')=\sum_{L'\in {\mathcal
L}'}(-1)^{|L'|}\iii_{L'}^{G'}\r_{L'}^{G'}(\pi').$$
Aussi,  
$$i({\bf JL}(\pi'))=\sum_{L\in {\mathcal
L}}(-1)^{|L|}\iii_L^G\r_L^G({\bf JL}(\pi'))=S_1+S_2$$ 
o\`u $S_1$ est la somme sur les sous-groupes
de Levi standard de $G$ qui se transf\`erent et $S_2$ sur les
sous-groupes de Levi standard de $G$ qui ne se transf\`erent
pas. La deuxi\`eme de ces sommes est nulle modulo $\ss$.

Les sous-groupes de Levi standard de $G$ qui se transf\`erent
correspondent biunivoquement aux sous-groupes de Levi standard de
$G'$.  Par
ailleurs, si $L$ et $L'$ se correspondent on a
$(-1)^{|L|}=(-1)^{n-r}(-1)^{|L'|}$. On a aussi, si $L$ et $L'$ se
correspondent, 
$$\jlh(\iii_{L'}^{G'}\r_{L'}^{G'}(\pi'))=\overline{\iii_L^G\r_L^G({\bf
JL}(\pi')) }$$
car $\jlh$ est un isomorphisme d'alg\`ebres de Hopf
(th. \ref{jlh}.c)). Le r\'esultat s'ensuit.\qed

Nous en d\'eduisons alors la
\begin{prop}\label{duales}
Si deux r.e.c.i. se correspondent par
la correspondance de Jacquet-Langlands, leurs duales se correspondent
au signe pr\`es. 
\end{prop}

Ce n'est pas, en g\'en\'eral, un cas particulier de
la conjecture. C'est n\'eanmoins ce qui importe pour les
correspondances globales \'etablies \`a la section suivante.

\section{Quelques r\'esultats globaux et leur cons\'equences locales}

\def\q{{\mathbb Q}}
\def\a{{\mathbb A}}
\def\ld{L^2(\omega,\gf\bc\ga)}
\def\bc{\backslash}
\def\lz{L^2_0(\omega,\gf\bc\ga)}
\def\gf{G'(F)}
\def\ga{G'(\a)}
\def\zf{Z(F)}

\def\za{Z(\a)}

\subsection{Un lemme}

Soit $F$ un corps local non-archim\'edien.
 Soit $G'$ une forme int\'erieure de $GL_n$ sur $F$. 
On dit qu'une repr\'esentation lisse irr\'eductible de $G'$ est
{\it elliptique} si son caract\`ere n'est pas identiquement 
nul sur l'ensemble des
\'el\'ements elliptiques r\'eguliers. Le lemme suivant affirme que
 $G'$ v\'erifie l'hypoth\`ese locale du cor. 7.3 de [Ar]. 

\begin{lemme}\label{lm1}
Soient $P$ un sous-groupe parabolique propre de $G'$, $L$ un
sous-groupe de Levi de $P$  et $\pi$ une
repr\'esentation unitaire de $L$. Alors l'induite parabolique
$ind_P^{G'}\pi$ n'a aucun sous-quotient elliptique. 
\end{lemme}

{\bf D\'emonstration.} Soit $\sigma$ une repr\'esentation lisse
irr\'eductible elliptique de $G'$. $\sigma$ peut s'\'ecrire dans
le groupe de Grothendieck de $G'$ comme une somme de repr\'esentations
obtenues par induction parabolique
\`a partir de repr\'esentations \eci\ ([Ta], prop. 2.1). 
Ces repr\'esentations
ont le m\^eme support cuspidal que $\sigma$. Comme $\sigma$ est
elliptique, ces repr\'esentations ne peuvent pas \^etre toutes 
induites \`a partir de sous-groupes paraboliques propres. Il y a donc
parmi elles au moins une (en fait exactement une, mais peu importe)
\eci\ de $G'$. Donc, le support cuspidal de $\sigma$
est forc\'ement un segment de Tadi$\check{\rm c}$.

Maintenant, si l'induite de $\pi$ est irr\'eductible, elle ne
peut pas \^etre elliptique (tout en \'etant son seul
sous-quotient) et le lemme est clairement vrai. 

Si l'induite de $\pi$ est r\'eductible, alors tout
sous-quotient de $ind_P^{G'}\pi$ en est une sous-repr\'esentation
parce que $\pi$ est unitaire. En appliquant la r\'eciprocit\'e de
Frobenius \`a chaque sous-repr\'esentation $\pi_i$, on trouve que
$\pi$ est un sous-quotient de $res_P^{G'}\pi_i$. Ainsi, si
$ind_P^{G'}\pi$ est r\'eductible, parmi les sous-quotients de
$res_P^{G'}ind_P^{G'}\pi_i$, $\pi$ appara\^\i t avec une multiplicit\'e
strictement sup\'erieure \`a 1 (*). Si on suppose par l'absurde que
l'induite de $\pi$ contient
une repr\'esentation elliptique $\sigma$, alors le support cuspidal de
$\pi$, 
\'etant le m\^eme que celui de $\sigma$,
doit \^etre un segment de Tadi$\check{\rm c}$ d'un sous-groupe de
Levi $L$, d'apr\`es ce qui pr\'ec\`ede. Dans la restriction 
\`a $L$ de l'induite de $\pi$ tous les sous-quotients apparaissent
donc avec multiplicit\'e 1 (**).  Il y a une contradiction entre
(*) et (**) qui prouve le lemme.\qed

\subsection{Sur la correspondance de Jacquet-Langlands globale}

Soient $F$ un corps global et $G'$ une forme
int\'erieure de 
$GL_n(F)$. Pour toute place $v$ de $F$ on
note $F_v$ le compl\'et\'e de $F$ en
$v$ et on pose $G_v=GL_n(F_v)$ et $G'_v=G'(F_v)$. Soient $Z$ le
centre de $G'$ et $\omega$ un caract\`ere unitaire de $\zf\bc\za$. On note
$R$ 
la repr\'esentation r\'eguli\`ere (translation \`a droite) 
de $\ga$ dans $\ld$, o\`u $\a$ d\'esigne
l'anneau des ad\`eles de $F$. On appelle {\it
s\'erie discr\`ete} toute sous-repr\'esentation irr\'eductible de
$R$. On note $R_{disc}$ la sous-repr\'esentation de $R$ engendr\'ee par la
partie discr\`ete. Si $\pi$ est une s\'erie discr\`ete de $G'$
et $v$ est une place de $F$,
on note $\pi_v$ la composante locale de $\pi$ \`a la place $v$.

Fixons une place finie $v$ de $F$. 
Rappels : soient $g\in G_v$ et $g'\in G'_v$ ;
on dit que 
$g$ et $g'$ {\it se correspondent} si leurs polyn\^omes
caract\'eristiques 
sont \'egaux, et sans racine multiple sur une cl\^oture
alg\'ebrique de $F_v$. Soit $Grot(G_v)$
(resp.$Grot(G'_v)$) le groupe de Grothendieck des repr\'esentations
lisses de longueur finie de $G_v$ (resp.$G'_v$).
Soient $\pi\in
Grot(G_v)$ et $\pi'\in Grot(G'_v)$. On dit que $\pi$ et $\pi'$ {\it
se correspondent} si l'on a la relation entre les caract\`eres~:
\begin{equation}\label{caracters1}
\chi_{\pi}(g)=(-1)^{n-r_v}\chi_{\pi'}(g')
\end{equation}
si $g$ et $g'$ se correspondent. 

Nous utilisons l'abr\'eviation {\it r.c.i.} pour repr\'esentation
de carr\'e int\'egrable. \`A la ssect. \ref{inv} nous avons \'etudi\'e
l'involution $i'$ d\'efinie sur $Grot(G'_v)$ par la formule de
dualit\'e 
de [Au]. \`A partir de maintenant
 nous dirons qu'une repr\'esentation irr\'eductible
$\pi$ de $G'_v$ est la duale d'une repr\'esentation irr\'eductible 
$\sigma$ si $i'(\pi)$ est \'egal \`a $\sigma$ ou $-\sigma$.
La dualit\'e sur $G_v$ obtenue de fa\c{c}on analogue est la
dualit\'e de Zelevinski. Rappelons que la duale d'une repr\'esentation
cuspidale est elle-m\^eme.
Nous utilisons l'abr\'eviation {\it d.r.c.i.} 
pour {\it repr\'esentation duale d'une repr\'esentation
de carr\'e int\'egrable}. Notons $Ca(G_v)$, $Ca(G'_v)$
(resp. $Du(G_v)$, $Du(G_v')$)  
l'ensemble des classes de r.c.i. (resp. 
d.r.c.i.) de $G_v$ et $G_v'$.
Il existe une unique bijection $C$ de $Ca(G_v)$ sur 
$Ca(G'_v)$ telle que tout $\pi\in Ca(G_v)$ corresponde 
\`a $C(\pi)$ (correspondance de Jacquet-Langlands).  
Il existe aussi une unique bijection $D$ (ne pas confondre avec
l'alg\`ebre \`a division des sections pr\'ec\'edentes et qui
n'appara\^\i tra plus) de $Du(G_v)$ sur 
$Du(G'_v)$ telle que tout $\pi\in Du(G_v)$ corresponde 
\`a $(-1)^{n-r}D(\pi)$
(prop. \ref{duales}).

Soient $\pi\in
Grot(G_v)\cup Grot(G'_v)$
et $\pi'\in Grot(G'_v)$. On dit que $\pi$ et $\pi'$ {\it
se
correspondent faiblement} s'il existe $\epsilon\in \{-1,1\}$ tel que 
les caract\`eres v\'erifient~:
\begin{equation}
\chi_{\pi}(g)=\epsilon\chi_{\pi'}(g')
\end{equation}
si $g\in G_v$ et $g'\in G'_v$ sont elliptiques r\'eguliers
et se correspondent.
\begin{lemme}\label{mobius}
Si $\pi$ est une repr\'esentation irr\'eductible elliptique de $G_v$
ou de $G'_v$, alors $\pi$ correspond faiblement \`a une r.e.c.i. de
$G'_v$.
\end{lemme}

{\bf D\'emonstration.} Suppsosns pour simplifier que $\pi$ est une
repr\'esentation de $G'_v$ (sinon on transf\`ere par ${\bf LJ}$ De
$G_v$ \`a $G'_v$).
On sait que $\pi$ s'\'ecrit comme une
combinaison lin\'eaire de
repr\'esentaions induites \`a partir de r.e.c.i. de \sgls\ de $G'_v$. 
Comme elle est elliptique, il y a au
moins un \sgls\ qui n'est pas propre. Comme toutes les r.e.c.i. qui
apparaissent ont le m\^eme support cuspidal, et comme au-dessus d'un
support cuspidal donn\'e il existe sur $G'_v$ au plus
une r.e.c.i. ([Ta]), on en d\'eduit qu'exactement l'un des \sgls\ n'est pas
propre. Comme les autres induites sont nulles sur l'ensemble des
\'el\'ements elliptiques r\'eguliers, on en d\'eduit que $\pi$
correspond faiblement \`a un multiple scalaire d'une
r.e.c.i. $\sigma$.
Le scalaire en question est le coefficient de $\sigma$ dans la somme
des induites de r.e.c.i. qui repr\'esente $\pi$. Ce coefficient est
plus ou moins un, par une formule d'inversion de M\"obius classique
(voir [Ze], prop 9.13).\qed

Soit $V$ l'ensemble des places de $F$ o\`u $G'$ est ramifi\'e. On suppose que
$V$ ne contient pas des places infinies. Rappelons qu'une repr\'esentation
$\pi$ de $G_v$ qu'elle est {\it $G'_v$-nulle}
 si le caract\`ere de $\pi$ est nul
sur l'ensemble des \'el\'ements de $G_v$ qui correspondent aux 
\'el\'ements de $G'_v$. Fixons deux places distinctes $v_1$ et $v_2$
de $F$.
On
note $A$ l'ensemble des s\'eries discr\`etes $\pi$
 de $G$ telles que, pour tout $v\in V$, $\pi_v$
n'est pas $G'_v$-nulle et la composante 
locale de $\pi$ en $v_1$ et $v_2$ est elliptique.

Soit $X_{\pi}$ l'ensemble des s\'eries discr\`etes 
de $G'$ telles que, pour toute 
place
$v\notin V\cup \{v_1,v_2\}$ on a $\pi_v\cong \pi'_v$. 

\`A partir de maintenant, sauf mention du contraire, le corps global
$F$ est {\it de caract\'eristique nulle}.

\begin{prop} \label{corresp1}

$X_{\pi}$ n'est pas vide.

\end{prop}

{\bf D\'emonstration.} Nous appliquons la d\'emarche 
classique de [JL], en
utilisant la formule des traces simple de [Ar], cor.7.5, avec des
fonctions \`a support dans les \'el\'ements elliptiques r\'eguliers
aux places $v_1$ et $v_2$. Cette formule s'applique \`a $G$ et $G'$
gr\^ace au lemme \ref{lm1}. On arrive \`a
\begin{equation}\label{formule}
\prod_{v\in V\cup \{v_1,v_2\}}\chi_{\pi_v}(g_v)=\sum_{\pi'\in 
X_{\pi}}
 m(\pi')\prod_{v\in V\cup \{v_1,v_2\}}\chi_{\pi'_v}(g'_v)
\end{equation}
($m(\pi')$ sont des multiplicit\'es non nulles) 
si, pour toute place $v\in V\cup \{v_1,v_2\}$,
$g_v\in G_v$ et $g'_v\in G'_v$ 
se correspondent et, de plus, $g_{v_1}$ et $g_{v_2}$
sont elliptiques r\'eguliers. Nous avons aussi utilis\'e
le th\'eor\`eme de multiplicit\'e un pour le spectre discret de $G$. 
Comme $\pi\in A$, on peut chosir les $g_v$ de telle fa\c{c}on que
$g_{v_1}$ et $g_{v_2}$ 
soient elliptiques r\'eguliers, que les autres $g_v$ correspondent \`a
des $g'_v\in G'_v$ et que le terme de gauche de la 
formule soit non nul. Donc le terme de droite n'est pas identiquement nul. 
Alors $X_{\pi}$ n'est pas vide.\qed

\begin{prop}\label{0}
Soit $\pi\in A$ cuspidale. Alors il existe $\pi'\in X_{\pi}$ telle que
$\pi'_v$ correspond \`a $\pi_v$ pour tout $v\in V\backslash\{v_1,v_2\}$,
 et $\pi'_v$ correspond faiblement \`a $\pi_v$ pour  
$v\in \{v_1,v_2\}$.
\end{prop}

{\bf D\'emonstration.} On utilise l'\'egalit\'e \ref{formule}. Par [Ba3], 
$X_{\pi}$ est fini.
\`A toute place 
$v\in V\backslash\{v_1,v_2\}$ la composante locale de $\pi$ est
g\'en\'erique ([Sh]).
Pour $v\in V\backslash\{v_1,v_2\}$ $\pi_v$ n'est pas $G'_v$-nulle par
hypoth\`ese. On applique la prop. \ref{generiques}. Le caract\`ere de
chaque $\pi_v$, $v\in V\backslash\{v_1,v_2\}$ peut \^etre ``remplac\'e''
dans la formule par le caract\`ere de la repr\'esentation $u'_v$
de $G'_v$
correspondante~:
$$\prod_{v\in V\bc \{v_1,v_2\}}u'_v(g'_v)
\prod_{v\in \{v_1,v_2\}}\chi_{\pi_v}(g_v)=
\epsilon \sum_{\pi'\in 
X_{\pi}}
 m(\pi')\prod_{v\in V\bc \{v_1,v_2\}}\chi_{\pi'_v}(g'_v)
\prod_{v\in \{v_1,v_2\}}\chi_{\pi'_v}(g'_v)$$
o\`u $\epsilon$ est un signe et, pour
$v\in \{v_1,v_2\}$, $g'_v$ correspond \`a $g_v$ elliptiques
r\'eguliers quelconques. Soit 
$X'_{\pi}$ l'ensemble de tous les $\pi'\in X_{\pi}$
tels que $\pi'_v$ est \'equivalent \`a $u'_v$ pour tout 
$v\in V\backslash\{v_1,v_2\}$. On peut alors \'ecrire~:
$$\prod_{v\in V\bc \{v_1,v_2\}}u'_v(g'_v)
\bigg( \prod_{v\in \{v_1,v_2\}}\chi_{\pi_v}(g_v)-\epsilon \sum_{\pi'\in 
X'_{\pi}}
 m(\pi')\prod_{v\in \{v_1,v_2\}}\chi_{\pi'_v}(g'_v)\bigg)=$$
$$=\epsilon \sum_{\pi'\in 
X_{\pi}\bc X'_{\pi}}
 m(\pi')\prod_{v\in V\bc \{v_1,v_2\}}\chi_{\pi'_v}(g'_v)
\prod_{v\in \{v_1,v_2\}}\chi_{\pi'_v}(g'_v)$$
pour tout $g_v$ et $g'_v$, $v\in \{v_1,v_2\}$, qui sont elliptiques
r\'eguliers et se correspondent. 
Par ind\'ependance lin\'eaire des caract\`eres sur
$$
\prod_{v\in
\{v_1,v_2\}}\chi_{\pi_v}(g_v)-\epsilon \sum_{\pi'\in  
X'_{\pi}}
 m(\pi')\prod_{v\in \{v_1,v_2\}}\chi_{\pi'_v}(g'_v)=0
$$ 
pour des tels
$g_v$
et $g'_v$.

Comme $\pi\in A$, il existe
 $g_{v_1}$ et $g_{v_2}$ elliptiques r\'eguliers
tels que 
$\chi_{\pi_{v_1}}(g_{v_1})$ et 
$\chi_{\pi_{v_2}}(g_{v_2})$ soient non nuls. 
Donc $X'_{\pi}$ est non vide.

Soit $X''_{\pi}$ l'ensemble des $\pi'\in X'_{\pi}$
tels que $\pi'_{v_1}$ et $\pi'_{v_2}$ sont elliptiques. 
Il est clair que la derni\`ere relation peut s'\'ecrire 
$\prod_{v\in V\backslash\{v_1,v_2\}}G'_v$, on obtient que 
$$\prod_{v\in
\{v_1,v_2\}}\chi_{\pi_v}(g_v)-\epsilon \sum_{\pi'\in  
X''_{\pi}}
 m(\pi')\prod_{v\in \{v_1,v_2\}}\chi_{\pi'_v}(g'_v)=0.$$

En utilisant le lemme \ref{mobius} et
l'orthogonalit\'e des caract\`eres sur $G'_v$, $v\in \{v_1,v_2\}$, 
on en d\'eduit
 que pour au moins une des $\pi'\in X''_{\pi}$, $\pi'_v$ correspond
faiblement \`a $\pi_v$, pour $v\in \{v_1,v_2\}$.\qed

Soit maintenant $\pi\in A$ dans le spectre r\'esiduel. Alors
$\pi$ provient d'un segment $[\rho,k]$ comme dans [MW], o\`u 
$k$ divise $n$ et $\rho$ est une repr\'esentation cuspidale
de $GL_{\frac{n}{k}}(F)$. Si $\rho_v$ est cuspidale, alors 
$\pi_v$ est une d.r.c.i..

\begin{prop}\label{1}
Supposons que $\rho_v$ est cuspidale pour toute $v\in V$.
Alors il existe $\pi'\in X_{\pi}$ telle que 

- pour toute
place $v\in V\backslash \{v_1,v_2\}$ on ait $\pi'_v=D(\pi_v)$
et 

- $\pi'_v$ correspond faiblement \`a $\pi_v$ pour $v\in \{v_1,v_2\}$.
\end{prop}

{\bf D\'emonstration.} La preuve est la m\^eme que pour la proposition 
pr\'ec\'edente, en utilisant l'application $D$ \`a la place de
l'application
$C$.\qed

\begin{prop}\label{22}
Supposons que $d$ divise $\frac{n}{k}$ et que pour toute 
place $v\in V$ $\rho_v$ est une r.c.i.. Alors il existe $\pi'\in
X_{\pi}$  
telle que 

- pour toute
place $v\in V\backslash \{v_1,v_2\}$ on ait 
$$\pi'_v=Lg(\nu^{-\frac{k-1}{2}}C(\rho_v)\otimes
\nu^{-\frac{k-3}{2}}C(\rho_v)) 
.... \otimes\nu^{\frac{k-1}{2}}C(\rho_v))$$
o\`u $\nu$ est la norme r\'eduite
et 

- $\pi'_v$ correspond faiblement \`a $\pi_v$ pour $v\in \{v_1,v_2\}$.
\end{prop}

{\bf D\'emonstration.} La preuve est la m\^eme que pour les propositions 
pr\'ec\'edentes. Seulement, quand on transf\`ere 
$\pi_v=Lg(\nu^{-\frac{k-1}{2}}\rho_v\otimes \nu^{-\frac{k-3}{2}}\rho_v
.... \otimes\nu^{\frac{k-1}{2}}\rho_v)$ on trouve un \'el\'ement de
$Grot(G'_v)$ dont l'auteur ne sait s'il est une repr\'esentation 
irr\'eductible ou pas (cas particulier qui semble assez difficile de
la conjecture pos\'ee). Si on \'ecrit cet \'el\'ement comme combinaison
lin\'eaire de repr\'esentations irr\'eductibles non-\'equivalentes,
alors nous avons prouv\'e, prop. \ref{corresp}, que 
$Lg(\nu^{-\frac{k-1}{2}}C(\rho_v)\otimes \nu^{-\frac{k-3}{2}}C(\rho_v))
.... \otimes\nu^{\frac{k-1}{2}}C(\rho_v))$ appara\^\i t dans cette somme
avec un coefficient non nul. Ceci suffit pour obtenir le r\'esultat.\qed
\ \\

{\bf Remarque.}
S'il s'av\`ere que la somme contient plus d'une repr\'esentation 
irr\'eductible i.e.
$Lg(\nu^{-\frac{k-1}{2}}\rho_v\otimes \nu^{-\frac{k-3}{2}}\rho_v
.... \otimes\nu^{\frac{k-1}{2}}\rho_v)$ ne correspond pas \`a
$Lg(\nu^{-\frac{k-1}{2}}C(\rho_v)\otimes \nu^{-\frac{k-3}{2}}C(\rho_v))
.... \otimes\nu^{\frac{k-1}{2}}C(\rho_v))$,
alors la conjecture de multiplicit\'e un 
est infirm\'ee pour les s\'eries discr\`etes de $G'$.
Ce probl\`eme semble \^etre assez important. Nous esperons 
qu'une correspondance de Jacquet-Langlands globale nous aide
\`a d\'emontrer un th\'eor\`eme de multiplicit\'e un pour $G'$ et \`a
classifier le sp\`ectre r\'esiduel de $G'$. Or, m\^eme si on applique
la formule des traces dans sa forme g\'en\'erale et on sait montrer 
l'\'egalit\'e des c\^ot\'es g\'eom\'etriques comme dans [AC], 
on se heurtera dans notre cas exactement
\`a ce probl\`eme (quel \'el\'ement de $Grot(G'_v)$ correspond \`a une
repr\'esentation irr\'eductible donn\'ee de $G_v$ ?) pour conclure sur
la multiplicit\'e ou sur la classification. D\'ej\`a sur le cas
particulier dans [Vi] on voit
qu'une correspondance de Jacquet-Langlands globale n'est pas induite
par les correspondances de Jacquet-Langlands locales, car aux places
o\`u $G'$ est ramifi\'e les r.c.i. de $G'_v$ 
correspondent tant\^ot  \`a des r.c.i. de $G_v$, tant\^ot \`a
des d.r.c.i..\\

Donnons un cas particulier, valable en toute caract\'eristique,
 utile peut-\^etre \`a  ceux qui ont de la
marge dans le
choix la forme int\'erieure ou la repr\'esentation $\pi$.
\begin{prop} Soit $\pi$ une s\'erie discr\`ete de $G$ telle que,
pour toute place $v\in V$, $\pi_v$ n'est pas $G'_v$-nulle. S'il existe
deux places $v_1$ et $v_2$ de $F$ telles que l'une des
deux conditions suivantes soit v\'erifi\'ee~:

(i) $\pi_{v_1}$ est cuspidale et $G'_{v_2}$ est le groupe des \'el\'ements
inversibles d'une alg\`ebre \`a divison,

(ii) $\pi_{v_1}$ et $\pi_{v_2}$ sont cuspidales,\\
alors il existe une unique s\'erie discr\`ete $\pi'$ de $G'$ telle que
$\pi_v$ et $\pi'_v$ se correspondent pour toute place $v$. Ceci est
vrai m\^eme si la caract\'eristique du corps global $F$ est non nulle. 
\end{prop}

{\bf D\'emonstration.} Comme $\pi_{v_1}$ est cuspidale, $\pi$ est
cuspidale.
On reprend la d\'emonstration de la prop. \ref{0},
en appliquant cette fois la
formule des traces simple de Kazhdan-Deligne, valable en toute
caract\'eristique, avec un coefficient non nul $f_{v_1}$
de $\pi_{v_1}$ \`a la
palce $v_1$ et une fonction $f_{v_2}$ \`a support dans l'ensemble des
\'el\'ements elliptiques r\'eguliers \`a la place $v_2$. 
D'une part,
$f_{v_1}$ n'annule pas la trace de $\pi_{v_1}$.
D'autre part, gr\^ace \`a l'hypoth\`ese sur la place $v_2$,
il existe un $f_{v_2}$ \`a support dans l'ensemble des \'el\'ements
elliptiques r\'eguliers de $G'_{v_2}$ qui n'annule pas la trace de
$\pi_{v_2}$. On obtient, comme dans la prop. \ref{0}, un $\pi'$ tel
que $\pi'_v$ correspond \`a $\pi_v$ pour tout $v\notin \{v_1,v_2\}$ et
$\pi'_{v_2}$ correspond faiblement \`a $\pi_{v_2}$. 
Toutefois, \`a cet instant, nous n'avons pas prouv\'e que $\pi'_{v_1}$
correspond faiblement \`a $\pi_{v_1}$, car la
fonction \`a la place $v_1$ a \'et\'e fix\'ee d\`es le d\'ebut. Mais on sait
quand m\^eme que la trace de $\pi'_{v_1}$ ne s'annule pas sur une fonction
qui correspond \`a $f_{v_1}$ (+). 
Soit $\sigma'_{v_1}$ la r.e.c.i. de $G'_{v_1}$ qui correspond \`a
$\pi_{v_1}$ par Jacquet-Langlands. 
Alors 

(++) $f_{v_1}$ correspond \`a un coefficient
$f'_{v_1}$ de $\sigma'_{v_1}$ et

(+++) comme $\pi_{v_1}$ est cuspidale,
$\sigma'_{v_1}$ est 
cuspidale (th. \ref{jlr}, b)).\\
De (+), (++) et (+++) on
obtient que $\pi'_{v_1}$ est \'equivalent \`a $\sigma'_{v_1}$.
Par cons\'equent, $\pi'_{v_1}$ est \'equivalent \`a $\sigma'_{v_1}$. 
Donc, $\pi_{v_1}$ correspond \`a $\pi'_{v_1}$. 

Regardons maintenant 
la place $v_2$ o\`u on sait que $\pi_{v_2}$ et $\pi'_{v_2}$
se correspondent faiblement. Dans l'hypoth\`ese (i), comme $\pi_{v_2}$
est \`a la fois elliptique et g\'en\'erique, c'est une r.e.c.i., et
comme $G'_{v_2}$ est le groupe des \'el\'ements inversibles d'une
alg\`ebre \`a division,  il est clair que $\pi_{v_2}$ et $\pi'_{v_2}$
doivent se correspondre par Jacquet-Langlands.
Dans l'hypoth\`ese (ii), on utilise le lemme \ref{mobius}~:

(a) $\pi'_{v_2}$ \'etant elliptique, elle correspond faiblement \`a une 
r.e.c.i. $\sigma'_{v_2}$ de $G'_{v_2}$,

(b) $\pi'_{v_2}$ et $\sigma'_{v_2}$ ont le m\^eme support cuspidal
(voir preuve du lemme \ref{mobius}).

Soit $\tau'_{v_2}$ l'image de $\pi_{v_2}$
par la correspondance de Jacquet-Langlands~; alors

(c) $\tau'_{v_2}$ est
alors cuspidale (th. \ref{jlr}, b)). 

Maintenant $\sigma'_{v_2}$ et $\tau'_{v_2}$ se correspond faiblement
puisque $\pi'_{v_2}$ et $\pi_{v_2}$ se correspond faiblement.
Comme $\sigma'_{v_2}$ et $\tau'_{v_2}$ sont des r.e.c.i., elles sont 
\'equivalentes (orthogonalit\'e des caract\`eres sur
$G'_{v_2}$). Par le (b) et (c), $\pi'_{v_2}$ est \'equivalente \`a 
$\tau'_{v_2}$.\qed

Le r\'esultat du corollaire suivant est d\'ej\`a connu
de tous. Le corps global $F$ est de caract\'eristique quelconque.
\begin{cor} Soit $S$ un ensemble fini de places finies de $F$. Pour
tout $v\in S$ on se donne une repr\'esentation de carr\'e int\'egrable
$\tau'_v$ de $G'_v$. Alors il existe une repr\'esentation automorphe
cuspidale $\pi'$ de $G'$ telle que, pour tout $v\in S$, $\pi'_v$ est
\'equivalente \`a $\tau'_v$. 
\end{cor}

{\bf D\'emonstration.} Quand $G'$ est $GL_n$, ce r\'esultat est, par
exemple, le lemme 6.5 de [AC]. Si $G'$ est une forme int\'erieure de
$GL_n$, posons, pour tout $v\in S$, $\tau_v=Ca^{-1}(\tau'_v)$.
On consid\`ere l'ensemble $\{\tau_v\}_{v\in S}\cup
\{\tau_{v_1},\tau_{v_2}\}$ o\`u $v_1$ et $v_2$ sont deux places qui ne
se trouvent pas dans $S$, et $\tau_{v_1}$ et $\tau_{v_2}$ sont des
repr\'esentations cuspidales unitaires de $G_{v_1}$ resp. $G_{v_2}$.
En appliquant le r\'esultat \`a $G=GL_n$, on trouve qu'il existe une
repr\'esentation cuspidale $\pi$ de $G$ telle que $\pi_v$ est
\'equivalente \`a $\tau_v$ pour tout $v\in S\cup \{v_1, v_2\}$. Si on
applique la proposition pr\'ec\'edente, cas (ii), \`a $\pi$, on trouve une
repr\'esentation cuspidale de $G'$ qui a les propri\'et\'es voulues.\qed

\subsection {Cons\'equences sur l'unitarisabilit\'e}

Nous donnons ici deux corollaires. Les r\'esultats dus \`a 
Tadi$\check{\rm c}$ ou 
conjectur\'es par Tadi$\check{\rm c}$ qui apparaissent dans cette
section se 
trouvent 
dans [Ta]. Le corps de base est ici de caract\'eristique nulle.

\begin{cor} Tout \'el\'ement de $Du(G'_v)$ est unitarisable.
\end{cor}

{\bf D\'emonstration.} On peut obtenir tout \'el\'ement de 
$Du(G_v)$ comme composante locale d'une repr\'esentation $\pi$ comme dans
la prop. \ref{1}, parce qu'on peut toujours trouver une repr\'esentation
cuspidale $\rho$ de $GL_{\frac{n}{k}}$ qui ait des composantes locales 
cuspidales unitaires donn\'ees aux places dans $V$.\qed
\ \\
\ \\
{\bf Remarque.} Ce corollaire
est un cas particulier de la conjecture $U1$ \`a la
fin de [Ta]
(ce que Tadi$\check{\rm c}$ y conjecture pour les r.c.i. nous l'avons
montr\'e ici pour les cuspidales unitaires).

\begin{cor} Si $G'_v=GL_r(A_v)$ ($A_v$ alg\`ebre \`a division)
 et l'entier $k$ divise $r$,
si $\delta$ est une r.c.i. de  $GL_{\frac{r}{k}}(A_v)$, alors
$$Lg(\nu^{-\frac{k-1}{2}}\delta\otimes \nu^{-\frac{k-3}{2}}\delta
.... \otimes\nu^{\frac{k-1}{2}}\delta)$$
est unitarisable.
\end{cor}

{\bf D\'emonstration.} On peut toujours trouver une repr\'esentation
cuspidale $\rho$ de $GL_{\frac{n}{k}}$ qui ait des composantes locales 
r.c.i. donn\'ees aux places dans $V$. La pr\'esente proposition est
une cons\'equence facile de ce fait, comme dans le corollaire pr\'ec\'edent,
mais en utilisant cette fois la prop \ref{22}.
\qed
\ \\
\ \\
{\bf Remarque.} C'est aussi un cas particulier des conjectures de 
Tadi$\check{\rm c}$ dans [Ta].
Prenons un exemple simple~:
$$Lg(\nu^{-2}\delta\otimes \nu^{-1}\delta
.... \otimes\nu^{2}\delta).$$
Supposons que $\nu_{\delta}=\nu^3$. 

Alors [Ta], prop. 4.3 implique
que $\nu^{-\frac{3}{2}}\delta\otimes \nu^{\frac{3}{2}}\delta$ est
r\'eductible et la conjecture $U1$ de Tadi$\check{\rm c}$ implique que 
$Lg(\nu^{-\frac{3}{2}}\delta\otimes \nu^{\frac{3}{2}}\delta)$ est 
unitarisable. 

La prop. 2.2 combin\'ee au lemme 2.5 (toutes les deux dans [Ta]) 
impliquent alors que
l'induite $\tau$ de
$\nu^{-\frac{1}{2}}Lg(\nu^{-\frac{3}{2}}\delta\otimes \nu^{\frac{3}{2}}\delta)
\otimes \nu^{\frac{1}{2}}Lg(\nu^{-\frac{3}{2}}\delta
\otimes \nu^{\frac{3}{2}}\delta)$ est
irr\'eductible. Elle est unitarisable par la conjecture $U2$ de
 Tadi$\check{\rm c}$.

Comme $\delta$ est unitaire aussi,
la premi\`ere conjecture $U0$ de Tadi$\check{\rm c}$ implique 
que l'induite de
$\tau\otimes\delta$ est irr\'eductible unitaire.

La prop 2.3 dans [Ta] implique que
l'induite de $\tau\otimes\delta$ n'est autre que 
$Lg(\nu^{-2}\delta\otimes \nu^{-1}\delta
.... \otimes\nu^{2}\delta).$
\ \\
\ \\

\section{Bibliographie}

[AC] J.Arthur, L.Clozel, {\it Simple Algebras, Base Change, and the
Advanced Theory of the Trace Formula},  Ann. of Math. Studies,
Princeton Univ. Press 120, (1989).

[Ar] J.Arthur, The invariant trace formula II. Global theory, J. of
the A.M.S. (1988), vol 1, no.3, 501-554.

[Au] A.-M. Aubert, Dualit\'e dans le groupe de Grothendieck 
de la cat\'egorie des repr\'esentations lisses de longeur 
finie d'un groupe r\'eductif $p$-adique, {\it Transactions A.M.S.}, 
Vol.347, No.6, 1995.

[Ba1] A.I.Badulescu, Orthogonalit\'e des caract\`eres pour
$GL_n$ sur un corps local de caract\'eristique non nulle, {\it
Manuscripa Math.} 101 (2000), 49-70.

[Ba2] A.I.Badulescu, Correspondance de Jacquet-Langlands en caract\'eristique
non nulle, pr\'epublication Univ. Poitiers no. 155, juin 2001, \`a 
para\^\i tre dans {\it Ann. Sci. de l'\'Ec. Norm. Sup.}

[Ba3] A.I.Badulescu, Un r\'esultat de finitude dans le spectre
automorphe pour les formes int\'erieures de $GL_n$ sur un corps global,
{\it Pr\'epublication}.

[BDK] J.Bernstein, P.Deligne, D.Kazhdan, Trace Paley-Wiener 
Theorem for reductive $p$-adic groups, {\it J. Analyse Math.} 
47 (1986), 180-192.

[Bo] Bourbaki, {\it Th\'eories spectrales, Chap.1-2}, Hermann, Paris.

[Ca] W.Casselman, Introduction to the theory of admissible
representations of reductive $p$-adic groups, preprint.

[Cl1] L.Clozel, Th\'eor\`eme d'Atiyah-Bott pour les vari\'et\'es
$p$-adiques et caract\`eres des groupes r\'eductifs, {\it M\'emoires
de la S.M.F.} no. 15, Nouvelle s\'erie, 1984, 39-64.

[Cl2] L.Clozel, Invariant harmonic analysis on the Schwarz
space of a reductive $p$-adic group, {\it Proc.Bowdoin Conf.1989,
Progress in Math.Vol.101}, Birkh\"auser, Boston, 1991, 101-102. 

[DKV] P.Deligne, D.Kazhdan, M.-F.Vign\'eras,
Repr\'esentations des alg\`ebres centrales simples $p$-adiques, {\it
Repr\'esentations des groupes r\'eductifs sur un corps local},
Hermann, Paris 1984. 

[vD] G. van Dijk, Computation of Certain Induced Characters of 
$p$-adic Groups, {\it Math. Ann.} 199, 1972, 229-240.

[FK] Y.Flicker, D.Kazhdan, Metaplectic correspondence, 
{\it Publ. Math. IHES} 64 (1987), 53-110. 

[H-CvD] Harish-Chandra,G. van Dijk, {\it Harmonic 
Analysis on Reductive $p$-adic Groups, L.N.M.}, Springer-Verlag, 1970.

[JL] H.Jacquet, R.P.Langlands, {\it Automorphic forms on
GL(2)}, Lecture Notes 114, Springer 1970.

[Le] B.Lemaire, Int\'egrales orbitales sur $GL(N)$ et corps 
locaux proches, {\it Ann.Inst. Fourier} 46 (1996), 1027-1056.

[MW] C. Moeglin, J.-L. Waldspurger, Le spectre r\'esiduel de $GL_n$,
Ann. Sci. \'Ecole Norm. Sup. t.22 (1989), 605-674.

[Sh] J.A.Shalika, The multiplicity one theorem for $GL_n$, {\it
Ann. of Math.} 100 (1974), 171-193.

[Ta] M.Tadi$\check{\rm c}$, Induced representations of $GL(n;A)$ for a 
$p$-adic division algebra $A$, {\it J. Reine angew. Math.} 405 (1990),48-77.

[Vi] M.-F. Vign\'eras, Correspondence between $GL_n$ and a division
algebra, {pr\'epublication}.

[Ze]  A.Zelevinski, Induced representations of reductive
$p$-adic groups II, {\it Ann. Sci. ENS} 13 (1980), 165-210.

\end{document}